\documentclass[12pt]{amsart}

\usepackage{color}
\usepackage{amsmath}
\usepackage{amsthm}
\usepackage{amssymb}
\usepackage{graphicx}
\usepackage{enumerate}
\usepackage{amsfonts}
\usepackage{mathrsfs}
\usepackage{parskip}
\usepackage{mathdots}
\usepackage{color}

\numberwithin{equation}{section}
\theoremstyle{plain}
\newtheorem{Proposition}[equation]{Proposition}

\newtheorem*{Corollary*}{Corollary}
\newtheorem{Theorem}[equation]{Theorem}
\newtheorem*{Theorem*}{Theorem}
\newtheorem{Lemma}[equation]{Lemma}
\theoremstyle{definition}
\newtheorem{Definition}[equation]{Definition}

\def\C{\mathbb{C}}
\def\R{\mathbb{R}}
\def\D{\mathbb{D}}
\def\T{\mathbb{T}}
\def\N{\mathbb{N}}
\def\Z{\mathbb{Z}}

\def\phi{\varphi}

\newcommand{\dist}{\operatorname{dist}}
\renewcommand{\ker}{\operatorname{Ker}}
\renewcommand{\dim}{\operatorname{dim}}

\renewcommand{\Re}{\operatorname{Re}}
\newcommand{\Lin}{\operatorname{Lin}}

\newcommand{\Pol}{\operatorname{Pol}}

\newcommand{\BMO}{\operatorname{BMO}}
\newcommand{\supp}{\operatorname{supp}}
\newcommand{\beqa}{\begin{eqnarray*}}
\newcommand{\eeqa}{\end{eqnarray*}}
\newcommand{\dst}{\displaystyle}
\newcommand{\lra}{\longrightarrow}
\newcommand{\lmto}{\longmapsto}
\newcommand{\Lra}{\Longrightarrow}
\newcommand{\Llra}{\Longleftrightarrow}

\title{Kernels of Toeplitz operators}

\author[Hartmann]{Andreas Hartmann}
\address{Institut de Math\'ematiques de Bordeaux, Universit\'e Bordeaux 1, 351 cours de la Lib\'eration 33405 Talence C\'edex, France}
\email{Andreas.Hartmann@math.u-bordeaux1.fr}

\author[Mitkovski]{Mishko Mitkovski} \address{Department of Mathematical Sciences\\
Clemson University\\ O-110 Martin Hall, Box 340975\\ Clemson, SC USA 29634} \email{mmitkov@clemson.edu}
\urladdr{http://people.clemson.edu/~mmitkov/}



\begin{document}

\begin{abstract}
Toeplitz operators are met in different fields of mathematics such as 
stochastic processes, signal theory, completeness problems, operator theory, etc. 
In applications, spectral and mapping properties are of
particular interest. In this survey we will focus on kernels of Toeplitz operators. This raises two 
questions. First, how can one decide whether such a kernel is non trivial? We will discuss in some details
the results starting with Makarov and Poltoratski in 2005 and their succeeding authors concerning this 
topic. In connection with these results we will also mention some intimately related applications to
completeness problems, spectral gap problems and P\'olya sequences. 
Second, if the kernel is non-trivial, what can be said about the structure of the kernel, and what kind of information on the Toeplitz operator
can be deduced from its kernel? In this connection we will review a certain number of
results starting with work by Hayashi, Hitt and Sarason in the late 80's on the extremal function.
\end{abstract}

\maketitle

\section{Introduction}
%


Toeplitz operators are natural generalizations of so-called Toeplitz matrices.
In the standard orthnormal basis of $\ell^2(\N)=\{a=(a_n)_{n\ge 0}:
\|a\|_2^2:=\sum_{n\ge 0}|a_n|^2<\infty\}$, a Toeplitz operator 
is represented by the infinite matrix
\begin{equation}\label{Toeplitzmatrix}
 T=\begin{pmatrix}
 u_0 & u_{-1} & u_{-2} & u_{-3} &\cdots\\
 u_1 &u_0 & u_{-1} & u_{-2}  &\ddots \\
 u_2 &u_1 & u_{0} & u_{-1}  &\ddots \\
 \vdots & \ddots & \ddots & \ddots & \vdots 
 \end{pmatrix},
\end{equation}
where $(u_n)_{n\in\Z}$ is a given sequence. 
If we identify $\ell^2(\N)$ with the usual Hardy space $H^2$ of functions
$f(e^{it})=\sum_{n\ge 0}a_ne^{int}$ with $(a_n)_n\in \ell^2(\N)$,
and if we associate (formally) $u$ with the Fourier series $\varphi(e^{it})=\sum_{n\in\Z}u_ne^{int}$ 
then
\[
 (Ta)_k=\widehat{\varphi f}(k), \quad k\ge 0,
\]
whenever that makes sense.
%
Hence, $Ta$ defines the sequence of Fourier coefficients of the function
\[
 P_+(\varphi f),
\]
where $P_+:L^2(\T)\lra L^2(\T)$, $P(\sum_{n\in \Z}a_ne^{int})=\sum_{n\ge 0}a_ne^{int}$ is the
so-called Riesz (or Szeg\H{o}) projection. 
We thus may define the Toeplitz operator with symbol $\varphi$ by
\[
 T_{\varphi}:H^2\lra H^2,\quad T_{\varphi}f=P_+(\varphi f),
\]
which is a continuous operator on $H^2$, and the associated matrix in the orthonormal basis
$\{z^n\}_{n\ge 0}$ of $H^2$ is given by \eqref{Toeplitzmatrix}.

Two special, but very important, Toeplitz operators are the shift operator $Sf=P_+(zf)=zf$ and its adjoint $S^*f(z)=P_+(\bar{z}f)$ 
--- the so-called backward shift operator. Any Toeplitz operator satisfies the following ``almost commuting'' relation 
\begin{equation}\label{toeplid}
 S^*T_{\varphi}S=T_{\varphi}, 
\end{equation}
and as a matter of fact this latter operator equality --- which is a kind of displacement condition
in the matrix reflecting the constance along diagonals --- characterizes Toeplitz operators.

A closely related class of operators which are also called Toeplitz operators can be defined in the so called  ``continuous case''. For most parts these theories are parallel to each other and studying one or the other case mostly depends on the person's taste. However, it should be noted that some problems are much more natural to be considered in the discrete setting while others are more suitable for the continuous one. For this reasons in this survey we will switch from one case to the other whenever appropriate. 

In the continuous case one starts with a function $u\in L^\infty(\R)$ and associates to it the convolution operator $T: L^2(\R_+)\to L^2(\R)$ defined, as usual, by $Tf=u*f$. If we want to consider $T$ as an operator from $L^2(\R_+)$ to itself it is natural, as in the ``discrete case" above, to project back from $L^2(\R)$ onto $L^2(\R_+)$ using the Riesz projection $P_+$ (now considered on the real line). The operator $T_u=P_+T: L^2(\R_+)\to L^2(\R_+)$ obtained in this way is again called a Toeplitz operator with symbol $u$.
It should be noted here that in this form $T_u$ is also known as a Wiener-Hopf operator. In order to recover the form of the Toeplitz operators discussed above it suffices to apply the Fourier transform. The ``Fourier transformed'' $T_u$ becomes an operator on the Hardy space in the upper half plane $H^2(\C_+)=\mathcal{F}L^2(\R_+)$, given by $T_\Phi f=P_+\Phi f$, where $\Phi=\hat{u}$ and $P_+:L^2(\R)\to H^2(\C_+)$ is the Riesz (or Szeg\H{o}) projection introduced above. In the continuous case the role of the shift operator $S$ is played by the translation semigroup $S(t): L^2(\R_+)\to L^2(\R_+), t\geq 0$, defined by $(S(t)f)(x)=f(x+t)$. In the Fourier domain this semigroup becomes $S(t):H^2(\C_+)\to H^2(\C_+), t\geq 0$, $S(t)f=e^{itx}f(x)$. Notice again that for each fixed $t\geq 0$ both $S(t)$ and $S^*(t)$ are in fact Toeplitz operators with very simple symbols $e^{itz}$ and $e^{-itz}$. Several results about Toeplitz operators can be viewed as perturbation results for these two special (but important) classes of Toeplitz operators. Finally, notice that the analog of \eqref{toeplid} in the continuous case is $S(t)^*T_\Phi S(t)=T_\Phi, t\geq 0$ which can again be taken as defining property for Toeplitz operators in the sense that any operator which satisfies this identity must be unitarily equivalent to a Toeplitz operator. This is the famous Nehari theorem which served as an inspiration for many important results in harmonic analysis and operator theory. 

Toeplitz operators, as introduced above, have numerous applications in signal theory,
stochastic processes, interpolation and sampling problems, etc., some of which will be overviewed
in Subsection \ref{AppToepKern} below and discussed in more details in Section~\ref{completeness}.
We refer to the monograph by B\"ottcher and Silbermann~\cite{BS} as a general reference on 
Toeplitz operators and their applications. 
In this survey we will be interested in a very specific topic in connection with Toplitz operators, 
namely the kernels of these Toeplitz operators --- so called --- Toeplitz kernels. We will try to give an overview of results about Toeplitz kernels --- some of which are very recent --- connected to the following two topics: 
\begin{enumerate}
\item Triviality of Toeplitz kernels: Given a Toeplitz operator $T_\varphi$, how to tell from the symbol $\varphi$ whether its kernel is trivial or non-trivial? 
\item Structure of Toeplitz kernels: In case of non-triviality we would like to understand the structure of a given Toeplitz kernel as a subset of the Hardy space $H^2$ or more generally $H^p$. 
\end{enumerate}
Besides being interesting objects to study intrinsically, as shown by Makarov and Poltoratski~\cite{MP1}, Toeplitz kernels are also interesting to study because of their intimate connections to numerous problems in complex and harmonic analysis, as well as mathematical physics. 
In the remaining part of this introduction we would like to give a brief overview of the topics that will be discussed in this paper. 

\subsection{Triviality of Toeplitz kernels} As just mentioned, the first natural property that one would like to understand about Toeplitz kernels is whether they are trivial or not. As far as we know, currently  there is no explicit, easily checkable, criterion for triviality of a Toeplitz kernel for general symbols $\varphi$. However, if we restrict the class of symbols $\varphi$, then the important results of Makarov and Poltoratski~\cite{MP2} provide an ``almost solution'' of the triviality problem. We will discuss these results in some detail in the last section of the paper. The class of symbols which can be treated by Makarov-Poltoratski techniques consists of unimodular symbols (unimodularity is no restriction of generality in the discussion of kernels of Toeplitz operators as we will see later) of the form $\overline{B}b$ for inner functions $B$ and $b$ having their zeros accumulating
only at one common point and which both have only one possible singular point mass at that point (this corresponds to meromorphic inner functions in the upper half plane). 
We will discuss more thoroughly a list of problems related with this triviality question in 
Subsection \ref{AppToepKern}.

\subsection{Structure of Toeplitz kernels}\label{StructKerT}
Once the question of non-trivia\-li\-ty is clarified, we now switch to the description of kernels of
Toeplitz operators, and the information that can be deduced from them. To motivate the results that follow, consider the special case of Toplitz operators $T_\varphi$ with symbols of the form $\varphi=\bar{\theta}$, where $\theta$ is some inner function. Then clearly, $\ker T_\varphi=K_\theta$, where $K_\theta=H^2\ominus \theta H^2$ is the model space generated by $\theta$. It is a well know result of Beurling that model spaces can be characterized exactly as those subspaces of $H^2$ which are invariant under the backward shift. This raises the following natural question: Does there exist some analogous characterization for Toeplitz kernels with general symbols? The first step towards the solution of this problem was made by Hitt \cite{Hi} who characterized the class of so called nearly invariant subspaces. From \eqref{toeplid} it is easy to check that every Toeplitz kernel is nearly invariant. However, not every nearly invariant subspace can be represented as a Toeplitz kernel. Recall that a nearly invariant subspace (with respect to $S^*$)
is a subspace $M$ of $H^2$ satisfying the following property
\[ 
 f\in M, f(0)=0\Lra S^*f\in M.
\]
Hitt showed that nearly invariant subspaces (in $H^2$) are of the form
\[
 M=gK^2_I,
\]
where $K^2_I$ 
is the model space generated by the inner function $I$ and $g$ is the extremal function of $M$, meaning that it maximizes the real part at
$0$ among all the functions in $M$ with norm one.
Extremal functions played a crucial role in the work of Hayashi \cite{Ha90} who used them to 
identify those nearly invariant subspaces which are exactly kernels of Toeplitz operators (see Section \ref{S2} for precise
definitions). Later, in \cite{SaKern} Sarason gave an approach to Hayashi's result using de Branges-Rovnyak spaces. 

As we will see below, one important consequence of this line of results is the fact that for Toeplitz kernels $\ker T_\varphi$ which are non-trivial we can always assume that $\varphi$ is a
unimodular function which can be represented as
\[
 \varphi=\frac{\overline{Ig}}{g}.
\]
This result represents the initial point in the Makarov-Poltoratski treatment of injectivity (see Lemma~\ref{MPlemma} below).

The extremal function $g$ appearing in Hitt's description of nearly invariant subspaces can also be used to decide whether the corresponding Toeplitz operator is onto \cite{HSS}. One can view the surjectivity problem for a Toeplitz operator $T_\varphi$ as a ``strong-injectivity'' problem for the adjoint Toeplitz operator $T_{\bar{\varphi}}$. Namely, it is well known that the surjectivity of any operator is equivalent to the left-invertibility of its adjoint. Left-invertibility, on the other hand, being equivalent to injectivity with closed range, can be naturally viewed as a type of ``strong injectivity''. So the Toeplitz operator $T_\varphi$ is surjective if and only if the Toeplitz operator $T_{\bar{\varphi}}$ is ``strongly injective'', i.e., the corresponding Toeplitz kernel is ``strongly trivial''. The properties of rigidity and exposed points (for which there is 
still no meaningful characterization available) as well as the Muckenhoupt $(A_2)$ condition are 
central notions here when considering the Hilbert-space situation $p=2$. The problem was also 
considered for the non-hilbertian situation where the extremal function does not have the same
nice properties \cite{HS} (see also Hitt's unpublished paper \cite{Hip}).
Related results for $p\neq 2$ were discussed by C\^amara and Partington and will be presented in
Section \ref{pneq2}. 

Bourgain factorization allowed Dyakonov to give another description of Toeplitz kernels on
Hardy spaces $H^p$: 
for every unimodular symbol $\varphi$ (unimodularity is no restriction of generality in the
discussion of kernels of Toeplitz operators as we will see later) the kernel $\ker_p T_{\varphi}$ 
of $T_{\varphi}$ considered now on $H^p$, is of the form
\begin{equation}\label{dyakbourgain}
 \ker_p T_{\varphi}=\frac{g}{b}(K^p_B\cap bH^p),
\end{equation}
where the tripel $(B,b,g)$ (in Dyakonov's terminology) 
consists of two Blaschke products $B$ and $b$, and a bounded analytic
function $g$ which is also boundedly invertible. Being in a more general situation here than
the hilbertian case, we shall define $K^p_I=H^p\cap I\overline{zH^p}$ (on $\T$)
for an inner function $I$.
Clearly, in that situation we can replace the symbol $\varphi$ by
\[
 b\overline{B}\frac{\overline{g}}{g}.
\]
We will give some interesting consequences of this
representation in Section \ref{bourgain}. Dyakonov is actually able to construct symbols such that
the kernel of the Toeplitz operator takes precise given dimensions in $H^p$ depending on the
values of $p$.

The representation \eqref{dyakbourgain} makes a connection with the completeness problems
discussed above. 
Completeness of a system of reproducing kernels is {\it via} duality equivalent to uniqueness
in the dual space.
Then, if $\ker_p T_{\phi}\neq \{0\}$, where $\phi=b\overline{B}{\overline{g}}/{g}$,
then this means that the zeros of $b$ form a zero sequence of $K^p_B$. 
Though in Dyakonov's result, $B$ is a Blaschke product, one can consider more general 
inner functions $I$.

Finally we would like to mention that a new connection between Toeplitz kernels and
multipliers between model space has been established in the recent preprint \cite{FHRmult}.
In this topic one is for instance interested in knowing whether the Toeplitz kernel, if non trivial,
contains bounded functions. This relates to Dyakonov's result above as well as to another result in the work of Makarov and Poltoratski (see \cite[Section 4]{MP1}).
They discuss criteria which ensure that 
if a given Toeplitz kernel is non-trivial in some $H^p$ (and more generally in the Smirnov class
$N^+$) then in a certain sense, increasing the size of $K_B$ in the Dyakonov representation
\eqref{dyakbourgain} ensures non-triviality for the smallest kernel in the $H^p$-chain, i.e. in
$H^{\infty}$. 



\subsection{Applications of Toeplitz kernels}\label{AppToepKern} Since the list of problems which can be translated into injectivity problems of Toeplitz operators is quite long, we will concentrate here our discussion to three specific problems. In the solution of each of these problems the Beurling-Malliavin densities play a very important role. As will become clear below the reason for this is the fact that the injectivity of many Toeplitz operators is closely dependent on these densities. 

The first problem concerns the geometry or basis properties of reproducing kernels in model spaces. It is a well known idea that one can use Toeplitz operators to study the basis properties of a the normalized reproducing kernels in model spaces. This idea goes back at least to the seminal paper by Hruschev, Nikolski, and Pavlov~\cite{HNP} who used the Toeplitz operator approach to finally settle the Riesz basis problem for non-harmonic complex exponentials. They also proved that most of the basis properties of sequences of normalized  reproducing kernels in model spaces can be described in terms of invertibility properties of an appropriate Toeplitz operator. Let us briefly recall their well known idea. Given a model space $K_{I}=H^2\ominus IH^2$, 
where $I$ is an inner function, and a sequence $\Lambda\subset \D$ (or in $\C_+$) the aim is to decide when a sequence of normalized reproducing kernels$\{k^I_\lambda\}_{\lambda\in\Lambda}$ is complete, a Riesz sequence, a Riesz basis, \dots in  $K_{I}$. It can easily be shown that under certain conditions on $\Lambda$ (see \cite[Chapter D4]{nik02}) this happens if and only if the Toeplitz operator $T_{{I}\overline{B_\Lambda}}$ has dense range, is injective with closed range (i.e.\ left invertible), invertible, \ldots, where $B_\Lambda$ denotes the Blaschke product with zeros set $\Lambda$. Indeed, if $P_I$ is the orthogonal projection from $H^2$ to $K_I$, then the basis properties translate into mapping properties of $P_I:K_{B_\Lambda}\lra K_I$ (under certain conditions on $\Lambda$). It is not difficult to check that those mapping properties are reflected in those of the Toeplitz operator $T_{{I}\overline{B_\Lambda}}$ (see \cite[Lemma D4.4.4]{nik02}). Note that $T_{{I}\overline{B_\Lambda}}$ having dense range is equivalent to $T_{\bar{I}B_{\Lambda}}$ being injective. Now, the completeness problem for non-harmonic complex exponentials~\cite{redheffer} can be viewed as a completeness problem for normalized reproducing kernels in a suitable model space. Therefore, it can be restated as an injectivity problem for a suitable Toeplitz operator. This very well know notoriously difficult problem inspired a great deal of results in mathematical analysis in the first half of the 20th century. Levinson gave in 1936 a sufficient condition for completeness in $L^p(I)$. After many unsuccessful attempts, the problem was finally settled by Beurling and Malliavin in 1967~\cite{BM2} relying heavily on their previous results~\cite{BM1}. There are several different proofs that appeared since~\cite{MNH, koosis1996leccons, koosisII, dBr} but none of them is much simpler than the original one. 
The achievement of Makarov and Poltoratski here was that they were able to adapt the deep ideas of Beurling and Malliavin to give
a metric characterization of injectivity for Toeplitz operators with very general symbols (much more general than symbols needed to solve the classical completeness problem for complex exponentials). As a consequence they provided a solution to the completeness problem for normalized reproducing kernels in a very general class of model spaces. 


The second problem we would like to mention here is the spectral gap problem which is related to the uncertainty principle in harmonic analysis~\cite{HJ}. One of the broadest formulations of the uncertainty principle says that a function (measure) and its Fourier transform cannot be simultaneously small. There are many mathematically precise versions of this heuristic principle depending on what kind of smallness one is interested in. In the classical gap problem one is interested in the gaps in the support of the measure and its Fourier transform. The heuristics says that these gaps cannot be simultaneously too big. In other words, if the support of the measure has gaps that are too big then then the support of its Fourier transform cannot have too large gaps. As shown by the second author and Poltoratski, for certain sets $X$, the gap $[0,a]$ where $\hat{\mu}$ vanishes for a measure $\mu$ supported on $X$ can be measured by Beurling-Malliavin densities. More precise results and extensions will be discussed in Section \ref{completeness}.

The third problem we want to present here is the P\'olya problem. Here we make a connection with the area of entire functions which received much interest in the past due to its intimate connections to the spectral theory of differential operators, signal processing, as well as analytic number theory. The P\'olya problem is a uniqueness problem and asks for a description of separated real sequences with the property that there is no non-constant entire function of exponential type zero (entire function that grows slower than exponentially in each direction) which is bounded on this sequence. Such sequences are called P\'olya sequences. This problem, which plays a crucial role in the resolution of several important problems in the recent years~\cite{ABB, BBB}, was resolved only recently by Poltoratski and the second author in~\cite{MiP1}. As in the problems discussed above Toeplitz kernels play an important role in the solution of this classical problem.  Again the characterization of P\'olya sequences is expressed in terms of Beurling-Malliavin densities.


\section{Basic definitions}\label{S2}
We denote by $H^p$, $0<p<\infty$, the classical Hardy space of analytic functions on the
unit disk $\D=\{z\in \C:|z|<1\}$, for which
\[
 \|f\|_p^p:=\sup_{0<r<1}\frac{1}{2\pi}\int_{-\pi}^{\pi}|f(re^{it})|^pdt<\infty,
\]
and $H^{\infty}$ are the bounded analytic functions on $\D$ equipped with the 
usual norm $\|f\|_{\infty}:=\sup_{z\in\D}|f(z)|$.
As usual, we identify functions $f$ in $H^p$ with their non-tangential boundary limits on 
$\T$ also denoted by $f$.
More precisely,
associating with $f\in H^p$ its non-tangential boundary function enables us to identify isometrically
\[
 H^p=\{f\in L^p(\T):\hat{f}(n)=0,n>0\}.
\]
In view of this observation, we will mostly not distinguish between $f$ and its boundary function, 
and when speaking about
$H^p$ in this paper we indifferently mean the Hardy space of holomorphic functions on the 
unit disk or its boundary values in $L^p(\T)$. 

Observe that every function $f\in H^p$, as an analytic function on $\D$, may be written as its
Taylor expansion
\[
 f(z)=\sum_{n\ge 0}a_nz^n.
\]
At least when $p\ge 1$, we can associate with $f$ on $\T$ its Fourier series, and as a matter
of fact, the coefficients $a_n$ turn out to be the Fourier coefficients of $f$ on $\T$. 
This identifies isometrically $H^2$ with $\ell^2$ by associating with $f=\sum_{n\ge 0}a_nz^n\in H^2$ 
the sequence $(a_n)_{n\ge 0}$ of its Fourier (or Tayler) coefficients.

\subsection{Projection and Toeplitz operators}
Given any function in $L^p(\T)$, $1<p<\infty$, which is expressed as its Fourier series
$f(e^{it})=\sum_{n\in \Z}a_ne^{int}$, we can associate its truncation to non-negative
Fourier coefficients:
\beqa
 P_+:L^p(\T)&\lra& L^p(\T),\\
    f(e^{it})=\sum_{n\in \Z}a_ne^{int}&\lmto&f(e^{it})=\sum_{n\ge 0}a_ne^{int}.
\eeqa
This projection is called the Riesz (or Szeg\H{o}) projection. As mentioned above,
for $p=2$, the norm in $L^2(\T)$ is expressed
as the $\ell^2$-norm of its Fourier coefficients, and it is plain that in this situation $P_+$ is
continous (orthogonal projection). It turns out that this projection is also continuous for
$1<p<\infty$. 

Another way of defining the Riesz projection is {\it via} reproducing kernels. 
For $\lambda\in\D$, set
\[
 k_{\lambda}(z)=\frac{1}{1-\overline{\lambda}z},\quad z\in \D.
\]
Then, for every $f\in H^p$ and for every $\lambda\in\D$,
\begin{equation}\label{repkern}
  f(\lambda)=\langle f,k_{\lambda}\rangle:=\frac{1}{2\pi}\int_{-\pi}^{\pi}
 f(e^{it})\overline{k_{\lambda}(e^{it})}dt,
\end{equation}
which actually makes sense for every $1\le p\le \infty$. Equation \eqref{repkern} explains
the wording ``reproducing kernel" for $k_{\lambda}$. We then have
\[
 P_+f(z)=\langle f,k_z\rangle,
\]
and this is well-defined for every $f\in H^p$, $1\le p\le \infty$, $z\in\D$. However, when
$f\in L^1(\T)$ (or $L^{\infty}(\T)$) then the function $P_+f$ may be in weak-$L^1(\T)$ (or in
$\BMO(\T)$). We will not discuss these issues here and concentrate mainly on $p=2$ or on
$1<p<\infty$.

Pick now a bounded function $\varphi\in L^{\infty}(\T)$, then the Toeplitz operator with 
$T_{\varphi}$ symbol $\varphi$ is defined by
\[
 T_{\varphi}f=P_+(\varphi f).
\]
When we consider $T_{\varphi}$ on $H^2$, then {\it via} the identification $H^2-\ell^2$ mentioned
earlier this is exactly the form we had discussed in the introduction and which
thus generalizes to $H^p$.
Note that we can also write 
\[
 T_{\varphi}f(z)=\langle \varphi f,k_z\rangle,
\]
which not only makes sense for $1<p<\infty$ but also when $p=1$ or $p=\infty$.
\\

A similar set-up can be given in the upper half plane. 

\subsection{Factorization and model spaces}
Let us turn back to Hardy spaces. It is well known that every function $f\in H^p$ maybe
decomposed in three factors
\[
 f=BSF,
\]
where $B$ is the Blaschke product carrying the zeros set $\Lambda$ of $f$ (counting multiplicities):
\[
 B(z)=\prod_{\lambda\in\Lambda}\frac{\overline{\lambda}}{\lambda}
 \frac{\lambda-z}{1-\overline{\lambda}z},\quad z\in D,
\]
(the normalization is chosen to make
the Blaschke factors positive in zero). 
Note that the product $B$ converges if and only if $\Lambda$ satisfies the Blaschke condition
\[
 \sum_{\lambda\in\Lambda}(1-|\lambda|^2)<\infty.
\]
The function $S$ is the singular inner function determined by
a positive singular measure $\mu$ on $\T$:
\[
 S(z)=\exp\left(-\int_{-\pi}^{\pi}\frac{e^{it}+z}{e^{it}-z}d\mu(e^{it})\right),\quad z\in\D.
\]
Finally, $F$ is an outer function which is uniquely determined by the modulus of $f$ on $\T$
\[
 F(z)=\exp\left(\int_{-\pi}^{\pi}\frac{e^{it}+z}{e^{it}-z}\log|f(e^{it})|dt)\right),\quad z\in\D.
\]
Observe that $|B|=|S|=1$ a.e.\ on $\T$. Functions with $|I|=1$ a.e.\ on $\T$ 
are called inner, and any inner function is
a product of Blaschke product and a singular inner function (possibly trivial). Inner functions 
play a central r\^ole in Beurling's characterization of shift invariant subspaces. Recall that
\[
 S:H^p\lra H^p,\quad Sf(z)=zf(z),
\]
is called the shift operator.
\begin{Theorem}[Beurling]
Let $M\subset H^2$ be a closed subspace, non-trivial ($M\neq H^2$ and $M\neq \{0\}$). Then $M$ is
shift invariant if and only if there is a (unique) inner function $I$ such that $M=IH^2$.
\end{Theorem}
We can also consider the adjoint shift $S^*=T_{\bar{z}}$ which on $\D$ is given by
\[
 S^*f(z)=\frac{f-f(0)}{z}, \quad z\in \D.
\]
Then passing to orthogonal complements we
deduce that a closed non-trivial subspace $M\subset H^2$ is $S^*$-invariant if and only if there
exists a (unique) inner function $I$ such that
\begin{equation}\label{model2}
 M=K^2_I:=(IH^2)^{\perp}=H^2\ominus IH^2.
\end{equation}
We will also use the notation $K_I=K^2_I$. These are the so-called model spaces which play a 
central r\^ole in the theory of function models for certain classes of Hilbert space contractions
(see e.g. \cite[Chapter C]{nik02}).
By a result of Douglas, Shapiro and Shields \cite{DSS}, these spaces can also be defined by existence of
so-called pseudo-continuations (which we won't discuss here). 
The definition \eqref{model2} is not adapted to a generalization to values $p\neq 2$. It is is easy to
check that denoting by $H^2_0=zH^2$ the $H^2$-functions vanishing at 0, then
\[
 K^2=H^2\cap I\overline{H^2_0},
\]
where the equality has to be interpreted on the circle via non-tangential limits
(the bar-sign here means complex conjugation). 
This definition immediately passes
to any value of $p$ (but these spaces behave strangely for $p<1$ as discussed by Aleksandrov, 
for which we refer to the monograph \cite{CR}):
\[
 K^p:=H^p\cap I\overline{H^p_0}.
\]
The spaces $K^p$ are widely studied spaces but are still far from being completely understood. 
They constitute a central building block for kernels of Toeplitz operators. 

Let us recall three simple examples of model spaces.

\begin{itemize}
\item When $I(z)=z^n$, then $K^p_I=K^2_I=Pol_{n-1}$ the analytic polynomials of degree
at most $n-1$;
\item When $I(z)=e^{a(z+1)/(z-1)}$, then $K_I$ is isometrically isomorphic to 
the Paley-Wiener space $PW_{a/2}$ of entire functions of exponential type $a/2$ which are
square integrable on $\R$.
\item When 
$I=B$ is a Blaschke product with simple zeros then 
$K_B^p=\Lin(k_{\lambda}:\lambda\in\Lambda)^{-}$ (closure in $H^p$, $1<p<\infty$). In the
special situation when the zeros of $B$ form an interpolating sequence, then every $K^p_B$-function
is an $\ell^p$-sum of normalized reproducing kernels $k_{\lambda}/\|k_{\lambda}\|_{p}$.
\end{itemize}


For $1<p<\infty$, the spaces $K^p_I$ 
can also be understood as ranges of the following projection (orthogonal for
$p=2)$
\[
 P_I=IP_-\overline{I},\quad P_-=Id-P_+,
\] 
where $I$ is an inner function, so that $K^p_I=P_IH^p$, and
\[
 k_{\lambda}^I(z)=P_Ik_{\lambda}(z)=\frac{1-\overline{I(\lambda)}I(z)}{1-\overline{\lambda}z}
\]
is the reproducing kernel.

\subsection{Elementary properties of Toeplitz operators}

\begin{equation}\label{norm}
\|T_{\phi}\|_{H^2\to H^2} = \|\phi\|_{\infty}
\end{equation}
(for $1<p<\infty$, $p\neq 2$, this equality is an equivalence);
\begin{equation}
 T_{\varphi}^*=T_{\overline{\varphi}}
\end{equation}
\begin{equation}\label{BS}
T_{\phi} T_{\psi} = T_{\phi \psi} \iff \mbox{$\phi \in \overline{H^{\infty}}$ or $\psi \in H^{\infty}$};
\end{equation}
\begin{equation}\label{innker}
\mbox{$u$ inner} \implies \ker{T_{\overline{u}}} = K_u;
\end{equation}
\begin{equation}\label{inv}
\phi \in H^{\infty} \implies T_{\overline{\phi}} K_u \subset K_u.
\end{equation}
Observe that \eqref{inv} implies that the model spaces enjoy the so-called $F$-property: 
If $w$ is an inner function then
\begin{equation}\label{77642bxtslk}
f \in K_u \; \; \mbox{and} \; \; f/w \in H^p \implies f/w \in K_u.
\end{equation}


\section{Invertibility and left-invertibility of Toeplitz operators}\label{invprop}

It is slightly surprising that the injectivity problem for Toeplitz operators received significantly less attention in the past in comparison with the invertibility problem. As part of the invertibility problem, the so called injectivity with closed range (left-invertibility) of Toeplitz operators has been thoroughly studied. This is why we begin by giving a brief overview of the most well-known results about invertibility and left-invertibility of Toeplitz operators which could be used as a guide in the study of the injectivity problem. 

Recall that an operator $A:X\to Y$ between
two Banach spaces is called left-invertible (bounded from below, injective with closed range) if there exists a constant $\epsilon>0$ such that for every
$x\in X$, $\|Ax\|_Y\ge \epsilon \|x\|$. Left invertibility of $T_{\varphi}$ is equivalent to
surjectivity of its adjoint $T_{\overline{\varphi}}$. The following well known result  gives a necessary condition for a Toeplitz to be left-invertible (we refer to \cite[p.65]{BS}, see also \cite[Theorem B4.2.7]{nik1})
\begin{Theorem}[Hartman-Wintner]
If $\varphi\in L^{\infty}$ but not invertible in $L^{\infty}$, then $T_{\varphi}$ is not
bounded below. 
\end{Theorem}
Hence, if $T_{\varphi}$ is left-invertible or onto, then we can assume that the symbol $\varphi$ is bounded away from zero, i.e., $0<\delta\le |\varphi|
\le M<\infty$. In that situation, $\log |\varphi|$ is obviously integrable and hence there is an outer
function  $g\in H^2$ such that $|g|=|\varphi|$ a.e.\ on $\T$. As a matter of fact $g$ is invertible in
$H^{\infty}$ and thus $T_g$ is invertible with inverse $T_{1/g}$. 
Hence the Toeplitz operator $T_{\varphi}$ behaves
like the Toeplitz operator $T_{\varphi/g}$ which has a unimodular symbol. 

Invertibility of Toeplitz operators has been studied in the 1960's by Devinatz
and Widom in the Hilbert space situation. Their result can be stated as follows
\begin{Theorem}
Let $\varphi\in L^{\infty}$. Then the following are equivalent
\begin{itemize}
\item $T_{\varphi}$ is invertible on $H^2$.
\item $\varphi^{-1}\in L^{\infty}$, and there exists $h\in H^2$ such that
\[
 \frac{\varphi}{|\varphi|}=\epsilon \frac{h}{|h|},
\]
and $|h|^2$ satisfies the Helson-Szeg\H{o} condition ($\epsilon$ is a uni-modular constant).
\end{itemize}
\end{Theorem}
Recall that a measurable function $w$ on $\T$ is said to satisfy the Helson-Szeg\H{o} condition
if 
\[
 w=e^{u+\tilde{v}},
\] 
where $u$ and $v$ are real valued bounded functions with $\|v\|_{\infty}<\pi/2$. As a 
matter of fact the Helson-Szeg\H{o} condition turns out to be equivalent to the
Muckenhoupt $(A_2)$ condition. We introduce here the general case of the $(A_p)$
condition. 
\begin{Definition}
Let $w$ be a positive measurable function on $\T$. 
If
\[
 \sup_{I}\left(\frac{1}{|I|}\int_I w^p(e^{it})dt\right)^{1/p}\cdot
  \left(\frac{1}{|I|}\int_I w^{-p'}(e^{it})dt\right)^{1/p'}<\infty,
\]
where $1/p+1/p'=1$ and $I$ runs through all the subarcs of $\T$, then we say that 
$w^p$ satisfies the $(A_p)$ condition.
\end{Definition}

Hunt, Muckenhoupt
and Wheeden have shown that the already introduced Riesz projection $P_+$
is continuous on $L^p(\T,w^pdm)$, $1<p<\infty$, if and only if $w^p\in (A_p)$.

It is possible to give more equivalent conditions for invertibility of Toeplitz operators. We just
mention here the relation with Hankel operators for which the Nehari theorem gives
$\|H_{\varphi}\|=\dist (\varphi, H^{\infty})$. Since Hankel and Toeplitz operators complete
each other: $\varphi f=T_{\varphi}f+H_{\varphi}f$ (the reader who does not know anything about
Hankel operators might take this last identity as a definition of a Hankel operator), and when
$\varphi$ is uni-modular we get an equivalent criterion for invertibility of Toeplitz operators
{\it via} Nehari's theorem. Indeed, it is necessary and sufficient that $\dist(\varphi,H^{\infty})<1$
and $\dist (\overline{\varphi},H^{\infty})<1$. We shall not use this connection too much
here.

Note also that in the special situation when $\varphi$ is unimodular then the condition
$\varphi^{-1}\in L^{\infty}$ is automatic and the condition on invertibility is essentially
that of writing $\varphi=\overline{h}/h$ where $|h|^2$ satisfies the Helson-Szeg\H{o} or
Muckenhoupt $(A_2)$ condition.

With the $(A_p)$ condition in mind, it is natural to ask whether an invertibility criterion can be found 
for Toeplitz operators on $H^p$. As expected, and shown by Rochberg \cite{Ro77}, invertibility
for Toeplitz operators on $H^p$ is characterized by the $(A_p)$-condition. In case $\varphi$
is unimodular this reads as: there is $g\in H^p$ such that $|g|^p\in (A_p)$ and
\[
 \varphi=\frac{\overline{g}}{g} \text{ a.e.\ on }\T.
\]
In this situation, it is possible to define explicitely an inverse at least on a suitable dense set by
\[
 gT_{1/\overline{g}}.
\]
Indeed, this operator makes sense on the dense set $g\Pol\subset H^p$, and by direct inspection it is
seen to be the inverse to $T_{\overline{g}/g}$ on this dense set (observe that the inverse is in 
general not a Toeplitz operator).


\section{Injectivity of Toeplitz operators and rigid functions}

In this section we start discussing some general elementary
injectivity properties of Toeplitz operators. 

The following well known lemma by Coburn provides a very interesting special property possessed by the class of Toeplitz operators (see e.g.\ \cite[Lemma B4.5.6]{nik02}).
\begin{Lemma}
Let $\varphi\in L^{\infty}$, then at least one of the kernels $\ker T_{\varphi}$ or $\ker
T_{\varphi}^*=\ker T_{\overline{\varphi}}$ has to be trivial.
\end{Lemma}

Notice that there exist Toeplitz operators $T_\varphi$ such that both $T_\varphi$ and $T_{\bar{\varphi}}$ are injective. Understanding such symbols is in the heart of the injectivity problem. 


For most applications one can concentrate on unimodular symbols. For such symbols we have 
\begin{equation}\label{inject1}
 g\text{ outer in }H^2 \Lra T_{g/\overline{g}}\text{ injective on }H^2. 
\end{equation}
This is easily seen by contradiction. Indeed, if $f\in \ker T_{g/\overline{g}}$ then
$(Ig/\overline{g})f=\overline{\psi}$ for some function $\psi\in H^2$, $\psi(0)=0$. Then
$H^1\ni Igf=\overline{g \psi}\in \overline{H^1_0}$ which is possible only when $gf=g\psi=0$, i.e.\
$f=0$. 

Note that when $I$ is inner, then $Ig/\overline{g}=(1+I)g/(\overline{(1+I)}g)$ also gives a
symbol of an injective Toeplitz operator. 

If we reverse the numerator and the denominator, then the situation changes crucially. Answering the injectivity problem for such symbols requires
the notion of rigidity which we will define now.

\begin{Definition}\label{defrigid}
A function $f\in H^1$ is rigid if it is uniquely determined in
$H^1$ up to a positive constant multiplier 
by its argument: for every $g\in H^1$, $f/g>0$ a.e.\ $\T$, implies $f=\lambda g$,
$\lambda\ge 0$.
\end{Definition}

This parallels in a way
the definition of outer functions which are determined by their moduli. Rigid functions
are necessarily outer, since the argument of any inner factor $I$ is for instance also given by
$(1+I)^2$ (which can be immediately seen from elementary geometric considerations).
It is known that $f$ is rigid in $H^1$ if and only if $f/\|f\|_2$ is an exposed point in the 
unit ball of $H^1$ (\cite{dLR}). A simple sufficient condition for rigidity is for instance $f\in H^1$ and
$1/f\in H^1$ but this condition is not necessary. 
In \cite{IN} the authors show that every polynomial with no zeros in $\D$ and simple zeros on $\T$
is rigid. Moreover for every such polynomial and $f\in H^{\infty}$ with $\|f\|_{\infty}\le 1$,
the function $p\circ f$ is also rigid.
However, no useful characterization for rigid functions is known.
It turns out, rigidity is the right notion for injectivity of $T_{\overline{g}/g}$:
\[
 g^2\text{ rigid in }H^1\Llra T_{\overline{g}/g}\text{ injective on }H^2.
\]
In fact this observation can also be used as an equivalent definition of rigid functions 
(see \cite[X-2]{Sa}). 

Rigid functions arise naturally in the frame of completely non-determi\-nistic gaussian processes. 
The idea is to relate the spectral measure $d\mu=wdm$ of the process, 
where $w$ is log-integrable so that $w=|f|^2$ for some outer function $f\in H^2$, and 
consider the map $T: L^2(wdm)\lra L^2$, $Th=fh$ which maps isometrically
the future of the process to $H^2_0$ and the past to $(f/\overline{f})\overline{H^2}$. The process
is completely non-deterministic if the future and the past have trivial intersection, which happens
exactly when $T_{\overline{f}/f}$ is injective. We refer to \cite{BJH} for more precise 
information on this topic.


\section{Kernels of Toeplitz operators and extremal functions}\label{kernext}

\subsection{Nearly invariant subspaces}
We start recalling a definition from the introduction.
\begin{Definition}
A closed non trivial subspace $M$ of $H^2$ is called nearly invariant, if
\[
 f\in M \text{ and }f(0)=0\ \Lra\ S^*f\in M.
\]
\end{Definition}

Any space of the form $GK^2_I$ with $G(0)\neq 0$ is automatically
nearly invariant as can be seen from the following simple argument: if $f=Gh\in GK^2_I$ and
$f(0)=0$, then necessarily $h(0)=0$. Hence
\begin{equation}\label{NearlyInv}
 S^*f=S^*(Gh)=\overline{z}Gh=G\overline{z}h=GS^*h\in GK^2_I.
\end{equation}
However, the space $GK^2_I$ is in general neither contained in $H^2$, nor closed even it it were
contained in $H^2$. For this we would need that
$|G|^2dm$ is a Carleson and a reverse 
Carleson measure for $K^2_I$ (the first condition guarantees that $GK_I^2$ is a subspace of 
$H^2$ and the second one that it is closed; we refer to \cite{Cohn}, \cite{TV}, \cite{BFGHR}).

In his work on invariant subspaces of $H^2$ on an annulus, Hitt \cite{Hi} described the nearly
invariant subspaces in a precise way (he also wrote an unpublished paper for the case $1<p<\infty$,
\cite{Hip}). 
Pick a nearly invariant subspace $M$ and associate the
extremal function which is the solution to the problem
\[
  \sup \{\Re g(0):g\in M, \|g\|_2\le 1\}.
\]
This problem has a unique solution. Indeed, the existence of a solution follows from
an argument based on normal families, and switching to the equivalent problem
\[
 \inf \{\|g\|_2:g(0)=1,g\in M\},
\]
which is a closed convex set in the strictly convex space $H^p$, $1<p<\infty$, we see that the solution
has to be unique.

Hitt observed that in the case $p=2$
it is possible to divide isometrically by the extremal function and that the
resulting space is $S^*$-invariant:
\begin{Theorem}[Hitt]
Let $M$ be a nontrivial nearly $S^*$-invariant subspace of $H^2$, with extremal function 
$g$. Then $M = gK_I^2$, where $I$ is some inner function vanishing at the origin. Furthermore, $g$ 
is an isometric divisor on $M$: $\|f/g\|_2 = \|f\|_2$ for all $f \in M$.
\end{Theorem}

We will discuss the situation when $p\neq 2$ in Section \ref{pneq2}.

The function $g$ appearing in this result is not necessarily outer. Observe that given any
nearly invariant subspace with extremal function $g$ and associated inner function $I$, then
picking another inner function $J$ with $J(0)\neq 0$, it follows from \eqref{NearlyInv} 
that $JgK_I$ is nearly invariant.

On the contrary,  when $I(z)=z$ and $g$ is any outer function, then the space $M=gK_I$ is
nearly invariant (as a matter of fact, the only function in $M$ vanishing at 0 is 
the zero function, the backward shift of which is trivially in $M$). Picking an inner function
$J$ with $J(0)=0$ leads to a space $JgK_I$ which is not nearly invariant. 

Sarason described in a precise way the extremal function $g$ using a de Branges-Rovnyak spaces
approach. This also allowed him to actually characterize all the isometric multipliers on a given
model space $K_I^2$. For this we need to associate two parameters with $g$. First it is clear that
the extremal function is of unit norm, $\|g\|=1$. Then the measure $|g(e^{it})|^2dt$ is a
probability measure. Now the function
\[
 \frac{1}{2\pi}\int_{-\pi}^{\pi}\frac{e^{it}+z}{e^{it}-z}|g(e^{it})|^2 dt
\]
has positive real part and hence is the Caley transform of a function $b$ in the unit ball of
$H^{\infty}$:
\begin{equation}\label{definb}
 \frac{1+b(z)}{1-b(z)}=G(z):=\frac{1}{2\pi}\int_{-\pi}^{\pi}\frac{e^{it}+z}{e^{it}-z}|g(e^{it})|^2 dt.
\end{equation}
This maybe reinterpreted as saying that $|g(e^{it})|^2dt$ is the Aleksandrov-Clark measure
associated with the function $b$. Observe that $|g(e^{it})|^2dt$ is absolutely continuous so that
$b$ is not an inner function (we won't go further here into de Branges-Rovnyak spaces). 
Now set \cite{SaNI}
\[
 a=\frac{2g}{G+1},
\]
and observe that
\[
 b=\frac{G-1}{G+1}.
\]
Then
\[
 |a|^2+|b|^2=\frac{4|g|^2+|G|^2-2\Re G +1}{|G+1]^2},
\]
and since $|G|=|g|$ a.e. on $\T$, we obtain $|a|^2+|b|^2=1$ a.e.\ $\T$. As a result, every
normalized $H^2$-function $g$ can be written as
\[
 g=\frac{a}{1-b},
\]
where $a$ and $b$ are bounded analytic functions and $|a|^2+|b|^2=1$ a.e.\ $\T$.

Sarason's achievement in this context was to observe that whenever $I$ divides $b$, then setting
$b=Ib_0$, the function
\[
 g=\frac{a}{1-Ib_0}
\]
multiplies isometrically on $K_I$. We have to assume $I(0)=0$ which ensures that the extremal function
is in the space as it should. Observe that if $g$ multiplies isometrically on $K_I$, where $I(0)=0$, then
it is automatically the extremal function as seen from the following equality
(we can assume that $\Re(g(0))>0$):
\beqa
\lefteqn{ \sup\{|gh(0)|:gh\in gK^2_I,\|gh\|=1\}}\\
  &&=g(0)\sup\{|h(0)|:h\in K^2_I,\|h\|=1\}=g(0)
\eeqa
(the last observation follows from the fact that $|h(0)|\le \|h\|\le 1$ and $1\in K_I$).
Let us state this observation as
a separate result \cite[Theorem 2]{SaNI}.

\begin{Theorem}[Sarason]
Let $I$ be any inner function, $I(0)=0$. Then for every couple $a$, $b$ of bounded holomorphic 
functions on $\T$ with $|a|^2+|b_0|^2=1$ the function 
\[
 g=\frac{a}{1-Ib}
\] 
multiplies isometrically on $K_I$ (and is thus an extremal function of $gK_I^2$.
\end{Theorem}

\subsection{Kernels of Toeplitz operators}
Let $\varphi\in L^{\infty}(\T)$, and consider the Toeplitz operator $T_{\varphi}$. Since a 
Toeplitz operator is also defined by the operator identity \eqref{toeplid} we easily see that
the kernel of a Toeplitz operator is nearly invariant. Let us formally check this observation already
mentioned earlier. Suppose $T_{\varphi}f=0$ and
$f(0)=0$, then
\[
 T_\varphi (S^*f)=(S^*T_{\varphi}S)(S^*f)=S^*T_{\varphi}f=0,
\]
since when $f(0)=0$ then $SS^*f=f$.

This introduces immediately the next question: is it possible to identify the kernels of
Toeplitz operators beyond arbitrary backward invariant subspaces? The answer to this question
was given by Hayashi \cite{Ha90}. It is here that we need again rigid functions. He
gave the following classification.

\begin{Theorem}[Hayashi]
Let $M$ be a nearly invariant subspace with associated extremal function $g$ and inner
function $I$ ($I(0)=0$). Then $M$ is the kernel of a Toeplitz operator if and only if the
function 
\[
 g_0=\frac{a}{1-b_0}
\]
is rigid. Here $g$ is given by the Sarason parameters as
\[
 g=\frac{a}{1-Ib_0}.
\]
\end{Theorem}

It is known that when $g_0$ is rigid so will be $g$ (see e.g. \cite[Corollary to Proposition 6]{SaKern}). 
However, the converse is false:
we refer to \cite[Section 3]{Ha90} for an example 
of a function $g=a/(1-Ib_0)$ such that $g^2$ is rigid but 
$g_0^2=(a/(1-b_0))^2$ is not, and the nearly invariant subspace $gK^2_I$ is not the kernel
of a Toeplitz operator.

Observe that in this situation $g$ is automatically outer (since rigid functions are outer as we
have seen earlier). Note however that it is not required that $b_0$ is outer, it can actually have inner
factors.

Moreover, with the functions $g$ and $I$ associated with the kernel of a Toeplitz operators
$T_{\varphi}$, $\varphi\in L^{\infty}$, it is easy to check that
\[
 \ker T_{\varphi}=gK_I=\ker T_{\overline{Ig}/g}.
\]
So, as already observed, 
when considering kernels of Toeplitz operators, we can always assume that $\varphi$ is a
unimodular function represented as
\[
 \varphi=\frac{\overline{Ig}}{g}.
\]
Note that we do not claim that an arbitrary Toeplitz operator $T_{\varphi}$ can be
represented by $T_{\overline{Ig}/g}$.


\subsection{Surjectivity}

As it turns out, the extremal function $g$ of the kernel of a Toeplitz operator contains
even more information about the operator itself. We have seen earlier that rigidity of $g_0^2$,
where $g_0=a/(1-b_0)$, guarantees that the nearly invariant subspace is the kernel of a
Toeplitz operator which can be chosen to be $T_{\overline{Ig}/g}$. If moreover $|g_0|^2$ 
satisfies the $(A_2)$ condition, then the Toeplitz operator $T_{\overline{Ig}/g}$ is onto. That is
the result in \cite{HSS}.

\begin{Theorem}[Hartmann-Sarason-Seip]
Let $\varphi$ be a unimodular symbol. If $\ker T_{\varphi}$ is non trivial and
$\ker T_{\varphi}=gK_I$, where $g=a/(1-Ib_0)$, then $T_{\varphi}$ is onto if and only 
if $|g_0|^2\in (A_2)$.
\end{Theorem}

This result also makes a connection with left invertibility since $T_{\varphi}$ is onto
if and only if $T_{\overline{\varphi}}$ is left-invertible. As discussed in the introduction,
left-invertibility of Toeplitz operators
plays a central r\^ole for describing Riesz sequences of reproducing kernels (or interpolating
sequences) in model spaces (see \cite[Section 4.4]{nik02}).

Note that rigidity of $g_0^2$ is equivalent to injectivity of $T_{\overline{g_0}/g_0}$ while
$|g_0|^2$ is equivalent to invertibility of $T_{\overline{g_0}/g_0}$, so the $(A_2)$ condition
is a stronger requirement. 

In order to construct a function $g_0$ such that $g_0^2$ is rigid but $|g_0^2|\notin (A_2)$, one
can use the fact that the $(A_2)$-condition is open, which means that when a weight satisfies
$(A_p)$ then it also satisfies $(A_r)$ for $r$ sufficiently close to $p$. It is easily seen that if
$g_0\in H^2$ and $1/g_0\in H^2$ but $1/g_0\notin H^{2+\epsilon}$ for any $\epsilon>0$,
then $g_0^2$ is rigid but does not define an $(A_2)$-weight.
For such an example we refer to \cite{HSS}.


\section{Bourgain factorization}\label{bourgain}

Another way of writing the kernel of a Toeplitz operator was explored by Dyakonov \cite{dyakkern}.
Using the Bourgain factorization he was able to prove the result below. We shall first introduce
some notation. When $T_{\varphi}$ is considered on $H^p$, then we denote by $\ker_p T_{\varphi}$ 
the corresponding kernel. By a triple we mean three functions $(B,b,g)$, where $B$ and $b$ are
Blaschke products and $g$ is an invertible function in $H^{\infty}$.

\begin{Theorem}[Dyakonov, 2000]
(i) For any $0\neq \varphi\in L^{\infty}$, there exists a triple $(B,b,g)$ such that
\begin{equation}\label{dyakkern2}
 \ker_pT_{\varphi}=\frac{g}{b}(K^p_B\cap bH^p), \quad 1\le p\le\infty.
\end{equation}
(ii) Conversely, given a triple $(B,b,g)$, one can find a $\varphi\in L^{\infty}$ for which 
\eqref{dyakkern2} holds true. In fact, it suffices to pick $\varphi=b\overline{B}\overline{g}/g$.
\end{Theorem}

The triple appearing in (i) of the above theorem is not unique.

The representation \eqref{dyakkern2} is interesting in that it is naturally related to the injectivity
problem and hence to the completeness problem in $K^p_I$.
Indeed, if $K^p_I\cap bH^p$ is non
trivial, then there is a function in $K^p_I$ vanishing on the zeros of $b$. 
We will address the completeness problem in Section \ref{completeness}.

Putting both representations together, i.e. the Hitt-Hayashi representation and the Dyakonov 
representation of the kernel of a Toeplitz operator, we obtain (in the case $p=2$):
\[
 GK_I=\ker T_{\overline{IG}/G}=\frac{g}{b}(K^2_B\cap bH^2),
\]
where $G=a/(1-Ib_0)$ is the extremal function of $\ker T_{\overline{IG}/G}$ and
$(B,b,g)$ is the associated Dyakonov triple. The triple $(I,1,G)$ found from Hitt's result is
in general not suitable for the Dyakonov description since $G$ is in general neither in $H^{\infty}$
nor $1/G$.

Starting from his result on kernels of Toeplitz operator, Dyakonov went on further constructing 
symbols $\varphi$ for which the dimension of the kernels $\ker_p T_{\varphi}$ varies
in a prescribed way depending on $p$. More precisely he has the following result.

\begin{Theorem}
Given exponents $1=p_0<p_1<\cdots<p_N=\infty$, and integers
$n_1>n_2>\cdots>n_N=0$, there exist Blaschke products $B$ and $b$ satisfying
\[
 \dim (K^p_B\cap bH^p)=n_j,\quad p\in [p_{j-1},p_j),
\]
for $j=1,\ldots,N$.
\end{Theorem}

The idea of the construction may be illustrated for $N=3$ by choosing
\[
 G(z)=\prod_{k=1}^{m_1}(z-\zeta_k)^{-1/p_1}\cdot \prod_{l=1}^{n_2}(z-\eta_l)^{-1/p_2},
\]
where $\zeta_k$ and $\eta_l$ are different points on the circle and $m_1=n_1-n_2$. Then 
$G\in H^p$, $p<p_1$ and $1/G\in H^{\infty}$ (implying that $G$ is rigid). It remains to
put $\varphi=\overline{z^{n_1}G}/G$. For more details we refer the reader to \cite{dyakkern}.


\section{More on the case $p\neq 2$}\label{pneq2}

Concerning the theory developped in Section \ref{kernext},
the situation is much less clear in the non-Hilbert situation. 
For instance the extremal function does no longer have the nice multiplier properties it had
for $p=2$. This will be discussed in Subsection \ref{extrfct} below. 

We will also discuss a more general notion of rigidity.
This can actually be defined for the case $q>0$ as it was defined for $q=1$ (see Definition
\ref{defrigid}, and replace $H^1$ by $H^q$),
i.e.\ a function $f\in H^q$ is rigid if it is uniquely determined in $H^q$ (up to
a positive constant multiplier) by its argument. Clearly this is equivalent to say that
for every function $g\in H^q$
\[
  g/f\ge 0 \text{ a.e. }\T\quad  \Lra \quad g=\lambda f\text{ for some }\lambda\ge 0.
\]

Concerning injectivity of specific Toeplitz operators, the following result can be shown 
by a similar argument as in the case $p=2$ (see also\cite[Theorem 5.4]{CP})
\begin{Theorem}\label{thmCP5.4}
The Toeplitz operator $T_{\overline{f}/f}$ is injective on $H^p$ if and only if $f^2$ is rigid
in $H^{p/2}$.
\end{Theorem}

With this notion of rigidity in mind, we will present some extensions of Hayashi's result
to $p\neq 2$ as discussed in \cite{CP} and \cite{CPJOT}
in Subsection \ref{rigid}. 
Those authors also 
observed that the notion of near
invariance can be replaced in the following way. Let $M$ be a closed subspace of $H^p$, then
$M$ is nearly invariant if
\[
 f\in M,\ f\overline{z}\in H^p\ \Lra \ f\in M.
\]
Since $f\overline{z}\notin H^p$ when $f(0)\neq 0$ the above observation is immediate. The idea is 
to replace the function $\overline{z}$ by more general functions, and call $M$ nearly $\eta$-invariant
if
\[
 f\in M,\ f\eta\in H^p\ \Lra \ f\in M.
\]
Note that by the $F$-property \eqref{77642bxtslk} kernels of Toeplitz operators are automatically
nearly $\overline{\theta}$-invariant whenever $\theta$ is inner. However, we will not 
discuss this matter in this survey.

We end the section with some results on minimal kernels.

\subsection{Extremal functions}\label{extrfct}
Extremal functions for general $p$ are defined in exactly the same way as in the case $p=2$, 
i.e.\ if $M=\ker_p T_{\varphi}$, then $G$ is the
unique solution to the extremal problem
\[
 \sup\{\Re G(0):g\in M, \|G\|_p\le 1\}.
\]
Concerning the extremal function of the kernel of a Toeplitz operator, as soon as $p\neq 2$ we
lose the nice isometric multiplier property. 
In general $G$ is even not an isomorphic multiplier as illustrated by
the following result (\cite{HS}), and which relies on a notion of variational identity for extremal 
problems valid for general $p$ (as can be found in \cite{Shap71}).

\begin{Theorem}[Hartmann-Seip]\label{thmHS}
Let $T_{\varphi}$ be a Toeplitz operator on $H^p$, $1 < p < \infty$, 
and $G$ the extremal function of $\ker_p (T_{\varphi} )$.
\begin{enumerate}
\item If $p \le 2$, then $GK_I^2 \subset \ker_p(T_{\varphi}) \subset GK_I^p$
and $\|f/G\|_p \le c_p\|f\|_p$ for every function $f \in \ker_p(T_{\varphi})$.
\item If $p \ge 2$, then $GK_I^p \subset \ker_p(T_{\varphi}) \subset GK_I^2$
and $\|f\|_p \le c_q\|f/G\|_p$ for every function $f \in \ker_p(T_{\varphi})$ 
($1/p + 1/q = 1$). 
\end{enumerate}
In general, none of these norm estimates can be reversed unless $p = 2$.
\end{Theorem}

We refer to \cite{HS} were explicit examples are constructed showing the failure of the
reverse inequalities in general.

Still, we can relate the function $G$ with the situation $p=2$. Since
$G\in \ker_p T_{\varphi}$, we
have $\varphi G=\overline{I}\overline{\psi}$, with $I$ an inner
function vanishing at 0 and $\psi$ an outer function in $H^p$. The function $I$ will be called
the associated inner function.

\begin{Theorem}\label{Hittfcts}
Let $T_{\varphi}$ be a Toeplitz operator on $H^p$, $1<p<\infty$, and $G$ the extremal
function of $\ker_p T_{\varphi}$ with associated inner function $I$. Then
$g=G^{p/2}$ is the extremal function of a nearly $S^*$-invariant
subspace of $H^2$ expressible as $gK^2_I$.
\end{Theorem}

Since $g$ is extremal for $gK^2_I$, from Sarason's description we know that
\[
 g=\frac{a}{1-Ib},
\]
where $a$ and $b$ are bounded analytic functions with $|a|^2+|b|^2=1$ a.e.\ on $\T$.
A natural question is to ask whether $gK^2_I$ is the kernel of a Toeplitz operator, or in other
words if $g_0^2$ is rigid in $H^1$ where $g_0=a/(1-b)$. An answer to this question, especially
when $K^p_I$ is not finite dimensional, would certainly
give some more insight in the structure of kernels of Toeplitz operators for $p\neq 2$.
There are some partial answers due to
Cam\^ara and Partington that will be discussed in the next subsection.


%
%
%

\subsection{Finite dimensional kernels of Toeplitz operators}\label{rigid}

In this subsection we discuss the special situation when $\ker_p T_{\varphi}$ is finite
dimensional. Recall from Theorem \ref{thmHS} that if $G$ is the extremal function of the kernel,
then $\ker_p T_{\varphi}$ is included between $GK^p_I$ and $GK^2_I$ (with the right order depending
on whether $p\ge 2$ or $p\le 2$). From that observation, the only way of $\ker_p T_{\varphi}$ to be 
finite dimensional is that $I$ is a finite Blaschke product. In this situation, we get in particular
that $K^p_I=K^2_I$ which is just a space of rational functions:
\[
 K^p_I=K^2_I= \left\{\frac{p(z)}{\prod_{j = 1}^{n} (1 - \overline{\lambda_j} z)}: p \in 
 \Pol_{n - 1}\right\},
\]
where $\lambda_j$ are the zeros of $I$ repeated with multiplicity, and $\Pol_{n-1}$ are
the analytic polynomials of order at most $n-1$. In particular
\[
 \ker_p T_{\varphi}=G K^2_I=GK^p_I=
 G  \left\{\frac{p(z)}{\prod_{j = 1}^{n} (1 - \overline{\lambda_j} z)}: p \in 
 \Pol_{n - 1}\right\}.
\]
Observing that
$h(z)=\prod_{j = 1}^{n} (1 - \overline{\lambda_j} z)$ is an outer function which is obviously
invertible in $H^{\infty}$, we can also write
\[
 \ker_p T_{\varphi}=\frac{G}{h}\Pol_{n-1}.
\]
See also \cite[Theorem 2.8]{CPJOT}. Note that replacing the denominator
$\prod_{j = 1}^{n} (1 - \overline{\lambda_j} z)$ by any other dominator 
$\tilde{h}:=\prod_{j = 1}^{n} (1 - \overline{\mu_j} z)$, where $\mu$ are zeros of some other 
finite Blaschke product with same degree as $I$ we get
\[
 \ker_p T_{\varphi}= G K^2_I=\tilde{G} K^2_B,
\]
where $\tilde{G}=G\tilde{h}/h$, and $\tilde{h}/h$ is invertible in $H^{\infty}$.

A more subtle question is to decide whether a given space of the form $G\Pol_{n-1}$ can be
the kernel of a Toeplitz operator. A partial answer to this question can be found in Dyakonov's
result Theorem \ref{dyakkern2} which states that when $G$ is boundedly invertible, then $GK_I$ is
the kernel of a Toeplitz operator (and in particular when $I(z)=z^n$ which gives $K_I=\Pol_{n-1}$). 
So the interesting case appears when $G$ is unbounded.
Let us consider the case $n=1$, then
we need that $\C G=\ker_pT_{\varphi}$. This immediately gives that $G$ is the extremal
function (thus necessarily outer) of the kernel and the associated inner function is $I(z)=z$. 
Hence
\[
 \ker_p T_{\varphi}=G K_I=\C G.
\]
Note that when $p=2$, then by Hayashi's and Sarason's result we know that
\[
 G(z)=\frac{a(z)}{1-zb(z)}
\]
where $|a|^2+|b|^2=1$ a.e.\ on $\T$, and $G_0^2=(a/(1-b))^2$ is rigid. Even in the more general 
finite dimensional 
situation (and for arbitrary $1<p<+\infty$) it turns actually out that the rigidity assumption is
not needed for $g_0^2$ but only for $g^2$ (this contrasts to the general situation
discussed for $p=2$).

In this connection we start citing \cite[Theorem 5.4]{CP} here for the disk.
\begin{Theorem}[Cam\^ara-Partington]\label{thmCPrigid}
Let $f\in H^p$ and $M=\C f$ the space generated by $f$. Then $M$ is the kernel of a 
Toeplitz operator if and only if $f^2$ is rigid in $H^{p/2}$.
\end{Theorem}

Observe that replacing $GK_{z}$ by $\tilde{G}K_{b_{\lambda}}$ as explained above 
conserves rigidity ($G$ and $\tilde{G}$ are simultaneously rigid or not).

The above result generalizes to finite dimensional spaces (see \cite[Theorem 3.4]{CPJOT}) which we again
cite for the disk.

\begin{Theorem}[Cam\^ara-Partington]
Let $M\subset H^p$ ($1<p<\infty$) be a finite dimensional subspace $\dim M=N<\infty$. Then
$M$ is the kernel of a Toeplitz operator if and only if $M=GK_{z^N}=G\Pol_{N-1}$ and $G^2$ is
rigid in $H^{p/2}$.
\end{Theorem}

Note that for a unimodular symbol $\varphi$
\cite[Lemma p.15]{MP1} states that $\dim\ker_p T_{\varphi}=n+1$ if and only if
$\dim\ker_p T_{b_i^n\varphi}=1$ where $b_i$ is the Blaschke factor of the upper half-plane
vanishing at $i$.

\subsection{Minimal kernels}

Another observation from the work of Cam\^ara and Partington concerns minimality of
kernels. For a given function $f\in H^p$, the minimal kernel associated with $f$ is the kernel
of a Toeplitz operator such that the kernel of any other Toeplitz operator which annihilates $f$ 
contains this minimal kernel. The following, in a sense natural, result holds (\cite[Theorem 5.1]{CP}).

\begin{Theorem}[Cam\^ara-Partington]
Given $f=Iu\in H^p$, $1<p<\infty$. Then the minimal kernel is given by
\[
 K_{min}(f)=\ker_pT_{\overline{Iu}/u}.
\]
\end{Theorem}
In other words, for every $\varphi\in L^{\infty}(\T)$, if $f\in \ker_pT_{\varphi}$, then
\[
 K_{min}(f)\subset \ker_pT_{\varphi}.
\]
Note that it is clear that $K_{min}(f)$ contains $f$. If $f^2$ were rigid, then we would have
$K_{min}(f)=\C f$ as seen in Theorem \ref{thmCPrigid}, and the above theorem would be trivial. So, the
interesting situation is when $f$ is not rigid. This still remains largely unexploited territory. We
refer to \cite{CP} where some examples are discussed and again some links with rigidity are 
established. Note that when $p=2$, and $u$ can be represented as $a/(1-Ib)$ with
$|a|^2+|b|^2=1$, and $u_0^2=(a/(1-b))^2$ is rigid, then $K_{min}(f)=uK^2_I$ and $u$ is
extremal for this kernel. In particular, the minimal kernel can have arbitrary dimension (it can be
finite or infinite dimensional). Theorem \cite[Theorem 5.2]{CP} gives a general criterion for
$K_{min}(f)$ to be finite dimensional.

\begin{Theorem}
Let $0\neq f=Iu\in H^p$. Then $K_{min}(f)$ is finite dimensional if and only if $I$ is a finite
Blaschke product and $\ker_pT_{\overline{u}/u}$ is finite dimensional.
\end{Theorem}

We refer to \cite[Theorem 3.7]{CPJOT} which gives a criterion for the latter kernel to be finite-dimensional in
terms of suitable factorizations of the symbol. Still, it is not clear how to obtain these
factorizations from arbitrary symbols. Note that for instance if $u=(1+\theta)u_0$, where
$\theta$ is inner so that $(1+\theta)u_0$ is outer, then
it is always possible to write
\[
 \frac{\overline{u}}{u}= \frac{\overline{(1+\theta)u_0}}{(1+\theta)u_0}=
 \frac{\overline{\theta u_0}}{u_0},
\]
which defines a symbol containing at least $u_0K_{\theta}$.


\section{Triviality of Toeplitz kernels and applications}\label{completeness}

In the preceding sections of this survey we have largely discussed the structure of (non trivial)
kernels of Toeplitz operators. It is now time to discuss the 
problem of deciding when the kernel of a Toeplitz operator is non trivial.
An answer to that question was proposed in the work by Makarov and Poltoratski who realized
that 
a resolution of several classical open problems in the area of complex and harmonic analysis can be reformulated in terms of injectivity of Toeplitz kernels. As already mentioned earlier, these include the completeness problem for a wide class of model spaces, but also the gap problem, the type problem, the determinacy problem, and several others. The goal of this section is to give an overview of these results and to explain the general ideas behind their proofs. The accent will be of course placed on the role of the Toeplitz kernels.   


This somewhat longer section will be organized in the following way. In the first part we will concentrate on the triviality problem for Toeplitz kernels. In the second part we will touch upon several of the above mentioned applications. 

As already mentioned in the introduction, a more natural setting for the results that will be discussed in this section is the continuous setting, so everything in this section will be done exclusively in the upper half-plane and the real line.  Let $T_U:H^2\to H^2$ be a Toeplitz operator with a unimodular symbol $U$. The most interesting case in applications is the one when the symbol is a quotient of two meromorphic inner functions, i.e., $U=\bar{\Theta}\Psi$ with $\Theta$, $\Psi$ meromorphic inner functions. Recall that an inner function is said to be meromorphic if it can be extended to a meromorphic function on $\C$. It is not hard to see that those are exactly the inner functions that can be represented as $S^a B_\Lambda$, where $S(z)=e^{iz}$ and $B_\Lambda(z)$ is the Blaschke product whose zero set $\Lambda$ has no accumulation points on the real line. We will keep using $S(z)$ (as Makarov and Poltoratski do) to denote the singular inner function $e^{iz}$. Every meromorphic inner function  $\Theta$ can be represented on the real line by  $\Theta=e^{i\theta}$, where $\theta:\R\to\R$ is some strictly increasing continuous branch of the argument of $\Theta$. Note that the fact that $\Theta$ is a meromorphic inner function implies that $\theta$ is real-analytic. Therefore, the symbols that are of most interest can be written as $U=e^{i\gamma}$, with $\gamma:\R\to\R$ being a real-analytic function of bounded variation (difference of two real-analytic  increasing functions). 

The fundamental problem that was solved by Makarov and Poltoratski was the injectivity problem for Toeplitz operators with this type of symbols. As we said earlier one would like to be able to say whether $T_U$ is injective or not just by looking at the argument $\gamma$ of the symbol $U$. Recall first the two extreme cases discussed in Subsection \ref{StructKerT}, namely, $U=\bar{\Theta}$ and $U=\Psi$. In the first case clearly $T_U=T_{\bar{\Theta}}$ is not injective, whereas in the second case $T_U=T_{\Psi}$ is injective. If we look at this in terms of the arguments it suggests that if the argument function $\gamma$ is decreasing, then we don't have injectivity, and if the argument function is increasing we have injectivity (actually more then injectivity). It turns out that this trivial guess is not that far from the truth, but making it precise requires a great deal of effort. The first step towards the goal is given by the following simple (but fundamental) lemma. Before stating this lemma we recall that each outer function $H\in H^2$ can be represented (on the real line) as $H=e^{h+i\tilde{h}}$, where $h=\log|H|$ and satisfies $h\in L^1(d\Pi)$ and $e^h\in L^2(\R)$. Here $d\Pi(t)=\frac{\dst dt}{\dst 1+t^2}$. Conversely, each function $h:\R\to \R$ satisfying the last two conditions defines an outer function $H$ determined on the real line by $H=e^{h+i\tilde{h}}$. We use $\tilde{h}$ in the above formulas (and throughout this section) to denote the Hilbert transform of $h\in L^1(d\Pi)$ defined by \[\tilde{h}(x)=\frac{1}{\pi} v.p. \int_\R \left(\frac{1}{x-t}+\frac{1}{1+t^2}\right)h(t)dt.\]

\begin{Lemma}[Makarov, Poltoratski]\label{MPlemma} A Toeplitz operator $T_U:H^2\to H^2$ with unimodular symbol $U$ is non-injective if and only if there exists an inner function $\Phi$ and an outer function $H\in H^2$ such that \[U=\bar{\Phi}\frac{\bar{H}}{H}.\]  Alternatively, in terms of arguments, $T_U$ is non-injective if and only if the real-analytic increasing argument $\gamma$ of $U=e^{i\gamma}$ can be represented in the form $\gamma=-\phi-\tilde{h}$, where $\phi$ is the argument of some meromorphic inner function and $h\in L^1(d\Pi)$, $e^h\in L^1(\R)$. 
\end{Lemma} 

Part of this lemma has already been discussed in \eqref{inject1}.

Roughly speaking this lemma says that a Toeplitz operator $T_U$ is not injective if $U$ is ``close to being equal'' to $\bar{\Phi}$ for some inner function $\Phi$ which is exactly one of the extreme cases discussed above. The phrase ``close to being equal'', that we used in the previous sentence is admittedly very imprecise. However, as we will try to explain below it is not as wrong as it looks. Indeed, one of the crucial parts of Makarov-Poltoratski's proof consists in showing that the term $\bar{H}/H$ can be discarded by paying only a very small penalty.  

The lemma above tells us that we need to find a way to tell when the argument of the symbol $\gamma$ can be represented in the form $-\phi-\tilde{h}$, where $\phi$ is the (increasing) argument of some meromorphic inner function and $h\in L^1(d\Pi)$, $e^h\in L^2(\R)$. To do this Makarov and Poltoratski first considered the following simpler problem: \emph{How to tell by looking at $\gamma$ whether it can be represented as $\gamma=d+\tilde{h}$, for some decreasing $d$ and $h\in L^1(d\Pi)$.} To solve this problem they used the Riesz sunrise construction and considered the portion of the real line $\Sigma(\gamma)\subset \R$ defined by \[\Sigma(\gamma)=\{x\in \R : \gamma(x)\neq \sup_{t\in [x,\infty)}\gamma(t)\}.\] 

It is not hard to show that this set is open and consequently can be represented as a countable disjoint union of open intervals $\Sigma(\gamma)=\cup I_n$. Furthermore, for each connected component $I_n=(a_n, b_n)$ we have $\gamma(a_n)=\gamma(b_n)$ and $\gamma(x)<\gamma(a_n)=\gamma(b_n)$ for all $x\in (a_n, b_n)$. As a side note we mention that this construction was used by F. Riesz to provide one of the first proofs that the maximal operator is of weak type $(1,1)$. It was later used as a basis for the famous Calderon-Zygmund decomposition (see e.g.~\cite{gar}). Let $\gamma^*:\R\to \R$ be a function defined by $\gamma^*(x):=\sup\{\gamma(t): t\in [x, \infty)\}$. Clearly $\gamma^*$ is a decreasing function which is constant on each of the intervals $I_n$. Moreover, the difference $\delta:=\gamma^*-\gamma$ 
is a non-negative function supported on $\Sigma(\gamma)$. This way we have a representation $\gamma=\gamma^*-\delta$ with $\gamma^*$ decreasing. It is left to examine when the difference function $\delta$ can be represented as $\tilde{h}$ for some $h\in L^1(d\Pi)$. It turns out that this is not always the case. However, under an appropriate assumption on $\Sigma(\gamma)$ (that will be given momentarily) it is never too far from being true. The idea is to compare $\delta$ to a similar but much simpler function $\beta$ which is also supported on $\Sigma(\gamma)$.  We define $\beta$ on each $I_n$ by the tent function $T_n$ whose graph is an isosceles right triangle with a base equal to $I_n$, i.e., $T_n(x)=\text{dist}(x, \R\setminus I_n)$. Define then $\beta:\R\to \R$ as a linear combination of the tent functions $\beta=\sum_n\epsilon_n T_n$ with a freedom to choose the coefficients later. It is easy to see that a sufficient condition for $\beta\in L^1(d\Pi)$ is the so called shortness condition: 

\begin{equation}\label{short}
\sum_n \frac{|I_n|^2}{1+\text{dist}(0, I_n)^2}<\infty,
\end{equation}

With this condition in place it is not hard to show 
that for any $\epsilon>0$ one can choose the coefficients $0<\epsilon_n<\epsilon$ so that 
\begin{enumerate}
\item $|\beta'|\lesssim \epsilon$, and
\item $\delta-\beta$ can be represented as a sum of atoms, i.e., it belongs to the real Hardy space $H^1_{Re}(d\Pi)$. 
\end{enumerate}
Since $H^1_{Re}(d\Pi)\subset \tilde{L^1}(d\Pi)$ this gives the following representation 
\[ \gamma(x)-\epsilon x=\gamma^*(x)+(\delta(x)-\beta(x))+(\beta(x)-\epsilon x).\] Obviously the last term can be written as a sum of a bounded and a decreasing function. Therefore, it can be absorbed by the first two terms. This way, under the shortness assumption, we obtain the desired representation of $\gamma(x)-\epsilon x$ for any $\epsilon>0$. This beautiful idea, to use atoms to show that a function is a Hilbert transform of an $L^1$-function, can be traced back to Beurling-Malliavin~\cite{BM2}. 

It is remarkable that the shortness condition above is almost necessary for $\gamma$ to have such a representation. Namely, it can be shown that if the shortness condition fails, i.e., if
\begin{equation}\label{long}
\sum_n \frac{|I_n|^2}{1+\text{dist}(0, I_n)^2}=\infty,
\end{equation}
 then for no $\epsilon>0$ the function $\gamma(x)+\epsilon x$ can be represented as a sum of a decreasing function $d$ and $\tilde{h}$, with $h\in L^1(d\Pi)$. 

To summarize, we have the following result as the first step towards the solution of the injectivity problem for Toeplitz kernels. It is sometimes called  ``little multiplier theorem''. 

\begin{Theorem}[Makarov-Poltoratski] Let $\gamma:\R\to\R$ be a real-analytic function of bounded variation (difference of two increasing functions). Assume also that $\gamma'$ is bounded from below. 
\begin{itemize}
\item[(i)] If $\Sigma(\gamma)$ does not satisfy the shortness condition~\ref{short}, then for every $\epsilon>0$ the function  $\gamma(x)+\epsilon x$ cannot be represented as a sum $d+\tilde{h}$ of a decreasing function $d$ and some $h\in L^1(d\Pi)$. 

\item[(ii)] If $\Sigma(\gamma)$ satisfies the shortness condition~\ref{short}, then for every $\epsilon>0$ the function  $\gamma(x)-\epsilon x$ can be represented as $d+\tilde{h}$ for some decreasing function $d$ and some $h\in L^1(d\Pi)$. 
\end{itemize}

\end{Theorem}

There is an interesting function-theoretic interpretation of this result. Namely, in analogy with usual Toeplitz kernels in $H^2$ one can also define Smirnov-Nevanlinna Toeplitz kernels $\ker^+T_U$ as sets of locally integrable functions $f\in N^+$ such that $U\bar{f}\in N^+$, where we use $N^+$ as usual to denote the Smirnov class on the upper half-plane. It can be shown that for a unimodular symbol $U=e^{i\gamma}$ the  Toepliz kernel $\ker^+T_U$ is trivial if and only if the argument  $\gamma$ can be represented as a sum $d+\tilde{h}$ of a decreasing function $d$ and some $h\in L^1(d\Pi)$. This way one can view the theorem above as a solution to the injectivity problem in the Smirnov-Nevanlinna case. 

To solve the injectivity problem in the Hardy case requires an even more ingenious idea which also goes back to the work of Beurling and Malliavin~\cite{BM1}.  It is sometimes called the ``big multiplier theorem''.  It shows that under some rather mild regularity assumptions on $\gamma$ the shortness condition above is enough (up to an arbitrary small $\epsilon$ gap) to determine whether a function $\gamma:\R\to\R$ can be represented in the form $\gamma=-\phi-\tilde{h}$, where $\phi$ is argument of some meromorphic inner function, $h\in L^1(d\Pi)$, and morover $e^h\in L^1(\R)$. Thus, in the view of  Lemma~\ref{MPlemma}, the shortness condition can be used as an almost necessary and sufficient condition to test injectivity of Toeplitz operators. To prove this one needs to address the following two problems: 
\begin{enumerate}
\item replace the decreasing function $d$ with a stronger requirement that $d=-\phi$, where $\phi$ is an argument of some meromorphic function, 
\item in addition to condition $h\in L^1(d\Pi)$ we would also need the condition $e^h\in L^2(\R)$. 
\end{enumerate}
The second problem is much harder and its solution lies much deeper. 

We will outline here only the main idea behind the proof skipping many of the technical difficulties. For details we refer to~\cite{MP2}. It should be noted that a similar idea was also used in~\cite{dBr}. The crucial part of this problem can be formulated in the following way: \emph{Given a non-negative real-analytic function $h\in L^1(d\Pi)$, with $\tilde{h}'\lesssim 1$, and $\epsilon>0$, find a function $m: \R\to \R$ such that $h\leq m$ and $\epsilon x-\tilde{m}(x)$ is essentially an increasing function.} Indeed, assume that we can show that such a function $m$ exists. Let $\gamma=d+\tilde{h}$, with $d$ decreasing and $h\in L^1(d\Pi)$. Then 
\[\gamma(x)-\epsilon x=d(x)-(\epsilon x-\tilde{m}(x))+(\tilde{h}(x)-\tilde{m}(x)).\] This implies that there exists an argument of a meromorphic inner function (not a finite Blaschke product) $\phi_1$ and a bounded function $u$ such that $\gamma=-\phi_1+u+\tilde{h}_1$, where $h_1=h-m\leq 0$. It is then not hard to show that any such function can be represented as $-\psi+\tilde{k}$, where $\psi$ as an argument of some inner function and $e^{k}\in L^p(\R)$ with $p<1$ (see \cite{MP2} for details). Finally, one can improve the condition $e^{k}\in L^p(\R), p<1$ to $e^{k_1}\in L^1(\R)$ by moving an appropriate part of $\psi$ to $\tilde{k}$. 
We now return to the main problem; to find a function $m\in L^1(d\Pi)$ such that $h\leq m$ and $\epsilon-\tilde{m}'(x)\geq o(1)$. The first step consists in showing that the a priori assumptions on $h$ imply that $h_0(x)=h(x)/|x|$ satisfies $\int h_0(x)\tilde{h}'_0(x)dx<\infty$. The second step is then to use the finiteness of the last integral to set up the following extremal problem: Minimize \[I(m_0)=\int h_0(x)\tilde{h}'_0(x)dx+\epsilon\int |x|h_0(x)d\Pi(x)\] over all non-negative functions $m_0\in L^1(d\Pi)$ satisfying $m_0\leq h_0$. Usual abstract arguments can be used to show that this extremal problem has a solution $m_0$. Finally, the way the extremal problem is set up allows one to show that $m(x)=|x|m_0(x)$ satisfies our desired conditions.

We can now state the ``almost characterization'' of the injectivity of Teoplitz operators in terms of the shortness condition.

\begin{Theorem}[Makarov-Poltoratski]\label{inject} Let $T_U: H^2\to H^2$ be a Toeplitz operator with a unimodular symbol $U=e^{i\gamma}$. Let $\gamma:\R\to\R$ be a real-analytic function such that $\gamma=\phi-\psi$ where $\phi$ is an argument of a meromorphic inner function and $\psi$ is an increasing function such that $|\psi'|\simeq 1$. 

\begin{itemize}
\item[(i)] If $\Sigma(\gamma)$ doesn't satisfy the shortness condition~\ref{short}, then for every $\epsilon>0$ the Toeplitz operator $T_V: H^2\to H^2$ with symbol $V=Ue^{i\epsilon \psi}$ is NOT injective.  

\item[(ii)] If $\Sigma(\gamma)$ satisfies the shortness condition~\ref{short}, then for every $\epsilon>0$ the Toeplitz operator $T_V: H^2\to H^2$ with symbol $V=Ue^{-i\epsilon \psi}$ is injective. 
\end{itemize}

\end{Theorem}

In fact, Makarov and Poltoratski proved a more general result which also includes Toeplitz operators for which the function $\psi$ appearing in the above theorem is allowed to satisfy $|\psi'(x)|\lesssim |x|^\kappa$ as $x\to\infty$. This can be viewed as an extension of the classical Beurling-Malliaving theorem.  

We now show how this result and its extensions were used in the recent solutions of several long-standing classical problems. More details as well as several other applications can be found in~\cite{Pol}.

\subsection{Completeness problem} 

We have already mentioned in the introduction that the completeness of a sequence of reproducing kernels $\{k_{\lambda}^I\}_{\lambda\in\Lambda}$ in a model space $K_I$ can be characterized by the injectivity condition of $T_{\bar{I}B_{\Lambda}}$. The classical problem for completeness of non-harmonic complex exponentials on $L^2[0,1]$ is equivalent to the completeness problem for the normalized reproducing kernels in the model space $K_S$ where $S(z)=e^{iz}$. Thus, in this classical case, Theorem~\ref{inject} can be reformulated in terms of the Beurling-Malliavin densities. Namely, 

\begin{Proposition} Let $\Lambda=\{\lambda_n\}\subset \R$ be a discrete sequence of real numbers. Let $\Theta=e^{i\theta}$ be some meromorphic inner function whose increasing argument $\theta$ satisfies $\theta(\lambda_n)=2n\pi$ for all $n$.  

\begin{itemize}
\item[(i)] If $|\theta'(x)|\simeq |x|^\kappa$, then $\sup\{a\geq 0: \ker T_{\bar{\Theta}S^a} \neq \{0\}\}$ is equal to $D_{BM}^-(\Lambda)$, the interior Beurling-Malliavin density of $\Lambda$. 

\item[(ii)] For general $\theta$ we have $\inf\{a\geq 0: \ker T_{\bar{S}^a\Theta} \neq \{0\}\}$ is equal to $D_{BM}^+(\Lambda)$, the exterior Beurling-Malliavin density of $\Lambda$. In other words, the radius of completeness  for the sequence of complex exponentials $\{e^{i\lambda_n x}\}$, i.e., the supremum of all $a>0$ for which the sequence is complete in $L^2[0,a]$ is equal to $D_{BM}^+(\Lambda)$.
\end{itemize}

\end{Proposition}

The same approach can be utilized to solve the completeness problems for other families of special functions which naturally show up as eigenfunctions of some classical singular Sturm-Liouville (Schr\"odinger) operators. Many of these operators are singular at the endpoints which prevents direct application of the classical Beurling-Malliavin result (which corresponds to case when the characteristic function $\Psi$ for the operator satisfies the condition $|\Psi'(x)|\lesssim 1$). In these instances the general form of Theorem~\ref{inject} which allows a polynomial growth of $|\Psi'(x)|$ needs to be used.

\subsection{Spectral gap and oscillation}
Recall from the introduction that in the classical gap problem one is interested in the gaps in the support of the measure and its Fourier transform. The heuristics says that these gaps cannot be simultaneously too big. In other words, if the support of the measure has gaps that are too big then then the support of its Fourier transform cannot have too large gaps. It is customary to call the gap in the support of the Fourier transform a \emph{spectral gap}. The most rudimentary form of this principle is the well known fact that a measure that is supported on a finite interval $[-a, a]$ (so it is zero on a big portion of the real line) cannot have a spectral gap of positive measure. Indeed, the Fourier transform of such a measure is an entire function which obviously cannot vanish on a set of positive measure on the real line. A more advanced version is the Riesz brothers theorem saying that a measure that is supported on $[0,\infty)$ (so still being zero on a big portion of $\R$) also cannot have any spectral gaps. The general problem --- the gap problem --- can be formulated in the following way. Given a closed set $X\subseteq \R$, determine the largest spectral gap that a measure supported on $X$ can have. More precisely, we would like to find a way to compute the so called \emph{gap characteristics} $G(X)$ of $X$ which is defined by \[ G(X):=\sup\{a\geq 0 : \exists \mu, \text{ supp}\mu\subseteq X, \hat{\mu}=0 \text{ on } [0,a]\}.\] 
There are several classical results that address this problem, especially the part when $G(X)=0$. The well known Beurling gap theorem~\cite{koosisII} says that if a non-zero measure $\mu$ is supported on a closed set $X$ whose complement $X^c$ is long in the sense of~\ref{long}, then $\hat{\mu}$ cannot vanish on any interval of positive length. Another well known result in this direction is due to Levinson~\cite{Lev} who showed that if the tail $M(x)=|\mu|(x,\infty)$ of a non-zero measure $\mu$ satisfies $\int \log M(x)d\Pi(x)=-\infty$, then again the measure $\mu$ cannot have a spectral gap of positive length. A more sophisticated results generalizing these two classical statements were proved by de~Branges~\cite{dBr}, and later Benedicks~\cite{Ben}. 

A first step towards the final solution of the gap problem was obtained by Poltoratski and the second author in~\cite{MiP1} where it was proved that for closed sets $X$ which are discrete and separated the gap characteristics of $X$ is equal to the interior Beurling-Malliavin density of $X$. 
\begin{Theorem}[Mitkovski-Poltoratski] If $X=\{x_n\}$ is a separated discrete set, i.e., $\inf_{m\neq n} |x_n-x_m|>0$, then $G(X)=D^-_{BM}(X)$.
\end{Theorem}  

This result was later extended by Poltoratski~\cite{PolGap} to arbitrary closed sets $X$. He proved that in the general case, besides the density of $X$, an additional subtle energy condition enters into play. The shortest way to formalize the energy condition is through the notion of  \emph{$d$-uniform sequences}. A real sequence $\Lambda=\{\lambda_n\}$ is \emph{$d$-uniform} if 1) it is regular with density $d$, i.e., there exists a sequence of disjoint intervals $\{I_n\}$ in $\R$ satisfying the shortness condition $\sum_n (|I_n|/(1+\text{dist}(0, I_n)))^2<\infty,$ and $|I_n|\to \infty$ as $n\to \pm\infty$, such that $$\#(\Lambda\cap I_n)-d|I_n|=o(|I_n|) \text{ as } |n|\to \infty.$$ and 2) it satisfies the following energy condition: there exists a short partition $\{I_n\}$ such that
 $$ \sum_n\frac{\#(\Lambda\cap I_n)^2\log_+|I_n|-E_{I_n}(dn_{\Lambda})}{1+\text{dist}(0,I_n)^2}<\infty,$$ where $E_I(\mu)=\iint_{I, I} \log|x-y|d\mu(x)d\mu(y),$ is the usual energy of the compactly supported measure $1_I(x)d\mu(x)$. The measure $dn_{\Lambda}$ entering in the energy condition above is the counting measure on $\Lambda$. The solution of the gap problem is given by the following theorem.
 
 \begin{Theorem}[Poltoratski \cite{PolGap}]\label{gap} For any closed set $X\subset\R$,
$$ G(X)=\pi\sup\{d: \exists\ \   d\text{-uniform sequence }\{\lambda_n\}\subset X\}.$$
\end{Theorem}

Very recently, the second author jointly with Poltoratski, refined these results even further~\cite{MiP2}, and obtained a generalization of the Beurling spectral gap theorem that strengthened this theorem by a factor of two. More precisely, in this paper, among other things, a metric description of the gap characteristic was obtained when the positive and the negative parts of the measure are supported in certain prescribed parts of the set. The gap characteristic of a pair of disjoint closed subsets $A$ and $B$ of $\R$ is defined by
\[G(A,B)=\sup \{ a>0 : \exists \mu\not\equiv 0, \hspace{0.1cm}  \text{supp} \mu^-\subseteq A, \hspace{0.1cm} \text{supp} \mu^+\subseteq B, \hspace{0.1cm}  \hat{\mu}\equiv 0 \text{ on } (-a, a)\}.\]
As in the case of $G(X)$, the description of $G(A,B)$ depends on two properties of $A$ and $B$: their density and their energy.  

\begin{Theorem}[Mitkovski-Poltoratski \cite{MiP2}]\label{oscillation} For any disjoint closed sets $A, B\subset\R$,
$$ G(A,B)=\pi\sup\{d: \exists\ \   d\text{-uniform sequence }\{\lambda_n\},\ \{\lambda_{2n}\}\subset A, \{\lambda_{2n+1}\}\subset B\}.$$
\end{Theorem}
As a simple consequence of this oscillation theorem one obtains a sharpening of the famous oscillation theorem of Eremenko and Novikov~\cite{EN} (their theorem solved an old problem of Grinevich from 1964 that was included in Arnold's list of problems).

\begin{Theorem}[Mitkovski-Poltoratski \cite{MiP2}]\label{Eremenko} If $\sigma$ is a nonzero signed measure with spectral gap $(-a, a)$  then there exists an $a/\pi$-uniform sequence $\{\lambda_n\}$ such that  $\sigma$ has at least one sign change in every $(\lambda_n, \lambda_{n+1})$.
\end{Theorem}

The crucial step towards the proof of both the gap theorem and the oscillation theorem above is the following result about annihilating measures for a very large class of de~Branges spaces. This is where Toeplitz kernels enter into play. 

\begin{Theorem}[Mitkovski-Poltoratski \cite{MiP2}]\label{annihilator} Let $B_E$ be a regular de~Branges space and let $\Phi(z)=\overline{E(\bar{z})}/E(z)$ be the corresponding meromorphic inner function. If $\mu$ is a  non-zero measure that annihilates $B_E$ then there exists a meromorphic inner function $\Theta$ such that $\{\Theta=1\}\subset \supp \mu^+$ for which the Smirnov-Nevanlinna kernel $\ker^+T_{\bar{\Phi}\Theta}$ is non-trivial. The same holds for the support of the negative part of $\mu$ as well.
\end{Theorem}

\subsection{P\'olya's problem} 

As discussed in the introduction, a P\'olya sequence is a separated real sequence with the property that there is no non-constant entire function of exponential type zero (entire function that grows slower than exponentially in each direction) which is bounded on this sequence. 
Historically, first results on P\'olya sequences were obtained in the work of Valiron~\cite{Val}, where it was proved that the set of integers $\Z$ is a P\'olya sequence. Later, in ignorance of the work of Valiron, this problem was
popularized by P\'olya, who posted it as an open problem. Subsequently many different proofs and generalizations were given (see for example section 21.2 of~\cite{Levin} or chapter 10 of~\cite{Boa} and references therein).

In his famous monograph~\cite{Lev} Levinson showed that if $|\lambda_n-n|\leq
p(n)$, where $p(t)$ satisfies $\int p(t)\log|t/p(t)|dt/(1+t^2) < \infty$ and some smoothness conditions, then
$\Lambda=\{ \lambda_n \}$ is a P\'olya sequence. In the same time for
each such $p(t)$ satisfying $\int{p(t)dt/(1+t^2)}=\infty$ he was
able to construct a sequence $\Lambda=\{ \lambda_n \}$ that is not
P\'olya sequence. As it often happens in problems from this area,
the construction took considerable effort (see \cite{Lev}, pp.
153-185). Closing the gap between Levinson's sufficient condition
and the counterexample remained an open problem for almost 25
years until de Branges~\cite{dBr} solved it assuming extra regularity of the sequence. These results have remained strongest for a very long time.

Jointly with A. Poltoratski~\cite{MiP1} the second author recently derived the following complete characterization of P\'olya
sequences.

\begin{Theorem}[Mitkovski-Poltoratski~\cite{MiP1}]\label{T1.2}
A separated real sequence $\Lambda=(\lambda_n)_{n\in\Z}$ is a P\'olya sequence if and only if its interior Beurling-Malliavin density $D_{BM}^-(\Lambda)$ is positive.
\end{Theorem}

Here are some details of how Toeplitz kernels enter in the solution of this problem. Let $F$ be an entire function of exponential type $0$ which is bounded (say by $M$) on $\Lambda$. Choose a meromorphic inner function $\Theta$ such that $\{\Theta=1\}=\Lambda$. This should be done carefully, but we won't go into technicalities here. The choice of $\Theta$ and the description of $ D^-_{BM}(\Lambda)$ in terms of Toeplitz kernels imply that $\ker T_{\bar{\Theta} S^c}$ is nontrivial for all $c>0$. As before $S$ denotes the innner function $S(z)=e^{iz}$. Pick some $h \in \ker T_{\bar{\Theta} S^c}$ with $L^2$-norm $1$. It can be shown that one has the following Clark-type representation for the function $hF$:
\begin{equation}\label{E2.2}
h(z)F(z)=\frac{1-\Theta(z)}{2\pi
i}\int\frac{F(t)h(t)}{t-z}d\sigma(t),
\end{equation}
where $\sigma$ is the Clark measure associated to the inner function $\Theta$. Furthermore for any $n\in \N$, $F^n$ is still an entire function of exponential type $0$ (bounded by $M^n$ on $\Lambda$) so that one has the same Clark-type representation for $hF^n$ as well. A simple estimate then provides a bound for $h(x)F(x)^n$ for all $x\in \R$. Taking the $n$-th root and using that $n$ is arbitrary help to get rid of $h$. This way we obtain that $F$ must be bounded on $\R$ which combined with the fact that $F$ is of zero type implies that $F$ is constant function.

\bibliographystyle{plain}

\bibliography{referencesTK1}

\def\cprime{$'$} \def\cprime{$'$} \def\cprime{$'$}
\begin{thebibliography}{10}

\bibitem{ABB}
E.~Abakumov, A.~Baranov, and Y.~Belov.
\newblock Radial limits and invariant subspac.
\newblock {\em Int. Math. Res. Notices}, page to appear.

\bibitem{BBB}
A.~Baranov, Y.~Belov, and A.~Borichev.
\newblock A restricted shift completeness problem.
\newblock {\em Journal of Functional Analysis}, 263(7):1887--1893, 2012.

\bibitem{Ben}
M.~Benedicks.
\newblock The support of functions and distributions with a spectral gap.
\newblock {\em Mathematica Scandinavica}, 55:285--309, 1984.

\bibitem{BM1}
A.~Beurling and P.~Malliavin.
\newblock On fourier transforms of measures with compact support.
\newblock {\em Acta Mathematica}, 107(3):291--309, 1962.

\bibitem{BM2}
A.~Beurling and P.~Malliavin.
\newblock On the closure of characters and the zeros of entire functions.
\newblock {\em Acta Mathematica}, 118(1):79--93, 1967.

\bibitem{BFGHR}
A.~Blandign{\`e}res, E.~Fricain, F.~Gaunard, A.~Hartmann, and W.~Ross.
\newblock Reverse {C}arleson embeddings for model spaces.
\newblock {\em J. Lond. Math. Soc. (2)}, 88(2):437--464, 2013.

\bibitem{BJH}
P.~Bloomfield, N.~P. Jewell, and E.~Hayashi.
\newblock Characterizations of completely nondeterministic stochastic
  processes.
\newblock {\em Pacific J. Math.}, 107(2):307--317, 1983.

\bibitem{Boa}
R.~P. Boas.
\newblock {\em Entire functions}, volume~5.
\newblock Academic Press New York, 1954.

\bibitem{BS}
A.~B{\"o}ttcher and B.~Silbermann.
\newblock {\em Analysis of {T}oeplitz operators}.
\newblock Springer Monographs in Mathematics. Springer-Verlag, Berlin, second
  edition, 2006.
\newblock Prepared jointly with Alexei Karlovich.

\bibitem{CPJOT}
M.~C. C{\^a}mara and J.~R. Partington.
\newblock Finite dimensional toeplitz kernels and nearly-invariant subspaces.
\newblock {\em J. Op. Theor.}

\bibitem{CP}
M.~C. C{\^a}mara and J.~R. Partington.
\newblock Near invariance and kernels of {T}oeplitz operators.
\newblock {\em J. Anal. Math.}, 124:235--260, 2014.

\bibitem{CR}
J.~Cima and W.~T. Ross.
\newblock {\em The backward shift on the {H}ardy space}, volume~79 of {\em
  Mathematical Surveys and Monographs}.
\newblock American Mathematical Society, Providence, RI, 2000.

\bibitem{Cohn}
W.~S. Cohn.
\newblock Carleson measures for functions orthogonal to invariant subspaces.
\newblock {\em Pacific J. Math.}, 103(2):347--364, 1982.

\bibitem{dBr}
L.~De~Branges.
\newblock Hilbert spaces of entire functions.
\newblock 1968.

\bibitem{dLR}
K.~de~Leeuw and W.~Rudin.
\newblock Extreme points and extremum problems in {$H_{1}$}.
\newblock {\em Pacific J. Math.}, 8:467--485, 1958.

\bibitem{DSS}
R.~G. Douglas, H.~S. Shapiro, and A.~L. Shields.
\newblock Cyclic vectors and invariant subspaces for the backward shift
  operator.
\newblock {\em Ann. Inst. Fourier (Grenoble)}, 20(fasc. 1):37--76, 1970.

\bibitem{dyakkern}
K.~M. Dyakonov.
\newblock Kernels of {T}oeplitz operators via {B}ourgain's factorization
  theorem.
\newblock {\em J. Funct. Anal.}, 170(1):93--106, 2000.

\bibitem{EN}
A.~Eremenko and D.~Novikov.
\newblock Oscillation of fourier integrals with a spectral gap.
\newblock {\em Journal de math{\'e}matiques pures et appliqu{\'e}es},
  83(3):313--365, 2004.

\bibitem{FHRmult}
E.~Fricain, A.~Hartmann, and W.~Ross.
\newblock Multipliers between model spaces.

\bibitem{gar}
D.~J.~H. Garling.
\newblock {\em Inequalities: a journey into linear analysis}.
\newblock Cambridge University Press, 2007.

\bibitem{HSS}
A.~Hartmann, D.~Sarason, and K.~Seip.
\newblock Surjective {T}oeplitz operators.
\newblock {\em Acta Sci. Math. (Szeged)}, 70(3-4):609--621, 2004.

\bibitem{HS}
A.~Hartmann and K.~Seip.
\newblock Extremal functions as divisors for kernels of toeplitz operators.
\newblock {\em J. Funct. Anal.}, 202(2):342--362, 2003.

\bibitem{HJ}
V.~Havin and B.~J{\"o}ricke.
\newblock {\em The uncertainty principle in harmonic analysis}.
\newblock Springer, 1994.

\bibitem{Ha90}
E.~Hayashi.
\newblock Classification of nearly invariant subspaces of the backward shift.
\newblock {\em Proc. Amer. Math. Soc.}, 110(2):441--448, 1990.

\bibitem{Hip}
D.~Hitt.
\newblock The structure of invariant subspaces of {$H^p$} of an annulus.
\newblock Preprint.

\bibitem{Hi}
D.~Hitt.
\newblock Invariant subspaces of {$H^2$} of an annulus.
\newblock {\em Pacific J. Math.}, 134(1):101--120, 1988.

\bibitem{HNP}
S.~V. Hru{\v{s}}{\v{c}}{\"e}v, N.~K. Nikol{\cprime}ski{\u\i}, and B.~S. Pavlov.
\newblock Unconditional bases of exponentials and of reproducing kernels.
\newblock In {\em Complex analysis and spectral theory ({L}eningrad,
  1979/1980)}, volume 864 of {\em Lecture Notes in Math.}, pages 214--335.
  Springer, Berlin-New York, 1981.

\bibitem{IN}
J.~Inoue and T.~Nakazi.
\newblock Polynomials of an inner function which are exposed points in {$H^1$}.
\newblock {\em Proc. Amer. Math. Soc.}, 100(3):454--456, 1987.

\bibitem{koosis1996leccons}
P.~Koosis.
\newblock {\em Le{\c{c}}ons sur le th{\'e}oreme de Beurling et Malliavin}.
\newblock Publications CRM, 1996.

\bibitem{koosisII}
P.~Koosis.
\newblock {\em The logarithmic integral}, volume~2.
\newblock Cambridge university press, 2009.

\bibitem{Levin}
B.~I. Levin.
\newblock {\em Distribution of zeros of entire functions}, volume~5.
\newblock American Mathematical Soc., 1964.

\bibitem{Lev}
N.~Levinson.
\newblock {\em Gap and density theorems}, volume~26.
\newblock American Mathematical Soc., 1940.

\bibitem{MP1}
N.~Makarov and A.~Poltoratski.
\newblock Meromorphic inner functions, {T}oeplitz kernels and the uncertainty
  principle.
\newblock In {\em Perspectives in analysis}, volume~27 of {\em Math. Phys.
  Stud.}, pages 185--252. Springer, Berlin, 2005.

\bibitem{MP2}
N.~Makarov and A.~Poltoratski.
\newblock Beurling-{M}alliavin theory for {T}oeplitz kernels.
\newblock {\em Invent. Math.}, 180(3):443--480, 2010.

\bibitem{MNH}
J.~Mashreghi, F.~Nazarov, and V.~Havin.
\newblock Beurling--malliavin multiplier theorem: The seventh proof.
\newblock {\em St. Petersburg Mathematical Journal}, 17(5):699--744, 2006.

\bibitem{MiP1}
M.~Mitkovski and A.~Poltoratski.
\newblock P{\'o}lya sequences, toeplitz kernels and gap theorems.
\newblock {\em Advances in Mathematics}, 224(3):1057--1070, 2010.

\bibitem{MiP2}
M.~Mitkovski and A.~Poltoratski.
\newblock On the determinacy problem for measures.
\newblock {\em Inventiones Mathematicae}, 202(3):1241--1267, 2015.

\bibitem{nik1}
N.~K. Nikolski.
\newblock {\em Operators, functions, and systems: an easy reading. {V}ol. 1},
  volume~92 of {\em Mathematical Surveys and Monographs}.
\newblock American Mathematical Society, Providence, RI, 2002.
\newblock Hardy, Hankel, and Toeplitz, Translated from the French by Andreas
  Hartmann.

\bibitem{nik02}
N.~K. Nikolski.
\newblock {\em Operators, functions, and systems: an easy reading. {V}ol. 2},
  volume~93 of {\em Mathematical Surveys and Monographs}.
\newblock American Mathematical Society, Providence, RI, 2002.
\newblock Model operators and systems, Translated from the French by Andreas
  Hartmann and revised by the author.

\bibitem{PolGap}
A.~Poltoratski.
\newblock Spectral gaps for sets and measures.
\newblock {\em Acta mathematica}, 208(1):151--209, 2012.

\bibitem{Pol}
A.~Poltoratski.
\newblock {\em Toeplitz Approach to Problems of the Uncertainty Principle},
  volume 121 of {\em CBMS}.
\newblock American Mathematical Soc., 2015.

\bibitem{redheffer}
R.~M. Redheffer.
\newblock Completeness of sets of complex exponentials.
\newblock {\em Advances in Mathematics}, 24(1):1--62, 1977.

\bibitem{Ro77}
R.~Rochberg.
\newblock Toeplitz operators on weighted {$H^{p}$} spaces.
\newblock {\em Indiana Univ. Math. J.}, 26(2):291--298, 1977.

\bibitem{SaNI}
D.~Sarason.
\newblock Nearly invariant subspaces of the backward shift.
\newblock In {\em Contributions to operator theory and its applications
  ({M}esa, {AZ}, 1987)}, volume~35 of {\em Oper. Theory Adv. Appl.}, pages
  481--493. Birkh\"auser, Basel, 1988.

\bibitem{SaKern}
D.~Sarason.
\newblock Kernels of {T}oeplitz operators.
\newblock In {\em Toeplitz operators and related topics ({S}anta {C}ruz, {CA},
  1992)}, volume~71 of {\em Oper. Theory Adv. Appl.}, pages 153--164.
  Birkh\"auser, Basel, 1994.

\bibitem{Sa}
D.~Sarason.
\newblock {\em Sub-{H}ardy {H}ilbert spaces in the unit disk}.
\newblock University of Arkansas Lecture Notes in the Mathematical Sciences,
  10. John Wiley \& Sons Inc., New York, 1994.
\newblock A Wiley-Interscience Publication.

\bibitem{Shap71}
H.~S. Shapiro.
\newblock {\em Topics in approximation theory}.
\newblock Springer-Verlag, Berlin-New York, 1971.
\newblock With appendices by Jan Boman and Torbj{\"o}rn Hedberg, Lecture Notes
  in Math., Vol. 187.

\bibitem{Val}
G.~Valiron.
\newblock Sur la formule d'interpolation de lagrange.
\newblock {\em Bull. Sci. Math}, 49(2):181--192, 1925.

\bibitem{TV}
A.~L. Vol{\cprime}berg and S.~R. Treil{\cprime}.
\newblock Embedding theorems for invariant subspaces of the inverse shift
  operator.
\newblock {\em Zap. Nauchn. Sem. Leningrad. Otdel. Mat. Inst. Steklov. (LOMI)},
  149(Issled. Linein. Teor. Funktsii. XV):38--51, 186--187, 1986.

\end{thebibliography}

\end{document}